\documentclass[11pt]{article}
\usepackage{geometry}                
\pdfoutput=1
\usepackage{graphicx}
\usepackage{hyperref}
\usepackage{amssymb}
\usepackage{amsmath}
\usepackage{amsthm}
\usepackage{amsfonts}
\usepackage{tikz}
\usetikzlibrary{calc, shapes, backgrounds}
\usepackage{caption}
\usepackage{subcaption}
\newcommand{\BE}{\begin{equation}}
\newcommand{\Av}{\operatorname{Av}}
\newcommand{\EE}{\end{equation}}

\usepackage[T1]{fontenc}
\title{Counting occurrences of patterns in permutations.}
\author{Andrew R Conway and Anthony J Guttmann}
\begin{document}

\maketitle
\begin{abstract}
We develop a new, powerful method for counting elements in a {\em multiset.} As a first application, we use this algorithm to study the number of occurrences of patterns
in a permutation. For patterns of length 3 there are two Wilf classes, and the general behaviour of these is reasonably well-known. We slightly extend some of the known results
in that case, 
and exhaustively study the case of patterns of length 4, about which there is little previous knowledge. For such patterns, there are seven Wilf classes, and based on extensive
enumerations and careful series analysis, we have conjectured the asymptotic behaviour for all classes. 

\end{abstract}




\section{Introduction}
Let $\pi$ be a permutation on $[n]$ and $\tau$ be a permutation on $[k].$ $\tau$ is said to occur as a {\em pattern} in $\pi$ if for some sub-sequence of $\pi$ of length $k$ all the elements of the sub-sequence occur in the same relative order as do the elements of $\pi.$ If the permutation $\tau$ does not occur in $\pi,$ then this is said to be a {\em pattern-avoiding permutation} or PAP.

If the permutation $\tau$ occurs $r$ times, it is said to be an {\em $r$-occurrence} of the pattern. Clearly, pattern-avoidance corresponds to the case $r=0.$

Let $s_n(\tau)$ denote the number of permutations of $[n]$ that avoid the pattern $\tau.$ Stanley and Wilf conjectured, and Marcus and Tardos \cite{MT04} subsequently proved, that for any pattern $\tau$ in $[k]$ the limit $\lim_{n \to \infty} s_n(\tau)^{\frac{1}{n}}=\lambda$ exists and is finite. This means that the number of PAPs grows exponentially with $n,$ whereas of course the number of permutations of $n$ grows factorially. 

For the more general problem of $r$-occurrences of a given pattern $\tau,$ a similar result holds, and the exponential growth rate $\lambda$ is independent of $r,$ proved by Mansour, Rastegar and Roitershtein \cite{MRR19}.

There are 6 possible permutations of length three, and the number of permutations of length $n$ avoiding any of these 6 patterns is given precisely by $C_n = \frac{1}{n+1} \binom{2n}{n}\sim \frac{4^n}{\sqrt{\pi n^3}},$ where $C_n$ denotes the $n^{th}$ Catalan number. That is to say, all 6 possible patterns have the same exponential growth-rate as PAPs. Alternatively expressed, there is only one Wilf class for length-3 PAPs.

For length-4 PAPs there are three Wilf classes. Typical representatives of the three classes are $1234,$ $1342$ and $1324.$ The generating function for the first two classes is known. 

In the first case, first proved by Gessel \cite{IG90}, the generating function is D-finite, satisfying a third-order linear, homogeneous ODE and $s_n(1234) \sim  \frac{81\sqrt{3} \cdot 9^n}{16\pi \cdot n^4}.$ 

In the second case, first proved by B\'ona \cite{MB97}, the generating function is algebraic, and $s_n(1342) \sim \frac{64\cdot 8^n}{243\sqrt{\pi} \cdot n^{5/2}}.$ 

The third case has not been solved, but extensive numerical work by Conway et al. \cite{CGZ18} suggests that $s_n(1324) \sim C \cdot \mu^n \cdot \mu_1^{\sqrt{n}} \cdot n^g,$ where $\mu \approx 11.598$ (and possibly $9+3\sqrt{3}/2$ exactly), $\mu_1 \approx 0.0400,$ and $g \approx -1.1.$ The appearance of the sub-dominant term $\mu_1^{\sqrt{n}}$ is referred to as a {\em stretched exponential} term. If present, that would suffice to prove that the generating function is not D-finite \cite{GP17}.

Thus these three Wilf classes have generating functions that are D-finite, algebraic, and (almost certainly) non-D-finite respectively.



In this article we  study the more general question of the behaviour of the generating function of $r$-occurrences of a given pattern in a permutation of length $n.$ The problem
of pattern-avoiding permutations thus corresponds to the case $r=0.$

 For this more general problem, it is known that, for patterns of length-3, there are  two Wilf classes, one corresponding to the two patterns $123$ and $321,$
and the other corresponding to the remaining four permutations of length 3. 

Surprisingly, the $r$-dependence of the two classes is quite different. Let $\psi_r(n)$ denote the number of permutations of length $n$ containing exactly $r$ occurrences of the nominated pattern. Then for the class 123 and 321 one has
\BE \label{class1}
\psi_r(n) = \frac{Q_{2r}(n)(2n-r)!}{(n+2r+1)!(n-r)!} \sim \frac{C_r\cdot 4^n}{n\sqrt{n\pi}}
\EE
where $Q_{2r}(n)$ is a polynomial with integer coefficients of degree $2r.$
For small values of $r,$ the amplitude coefficient $C_r$ appears to increase exponentially, growing seemingly like $\lambda^r,$ where the growth constant
$\lambda \approx 2.67.$  That is to say, the asymptotics remain unchanged, and only the amplitude, or premultiplying constant changes with $r.$ However this exponential increase in the amplitude cannot continue indefinitely, and indeed must decline toward zero as $r$ becomes large, as we explain below.

For the second class, corresponding to patterns 132, 231, 213, 312, Mansour and Vainshtein \cite{MV02} show that
\BE \label{class2}
\psi_r(n)=\frac{Q_r(n)(2n-3r)!}{n!r!(n-r-2)!} \sim \frac{4^{n-3r/2} \cdot n^{r-3/2}}{r!\sqrt{\pi}}
\EE
where $Q_r(n)$ is a polynomial of degree $3(r-1),$ whose leading-order term is precisely $n^{3(r-1)}.$ Here the asymptotics are quite different. The amplitude changes in a very regular way, but the subdominant power-law exponent increases by 1 as $r$ increases by 1.


Another noteworthy property of these $r$-occurrences is that, while the counting sequences $\psi_r(n)$ are known or conjectured to be Stieltjes moment sequences for the case $r=0,$ as discussed in Bostan et al. \cite{BEGM20}, this is not the case for $r > 0,$ for all the cases we have studied. We will report further on this in a forthcoming paper.

Aspects of this problem have been previously studied by several authors. For the class 132, Noonan and Zeilberger \cite{NZ96} conjectured the result for $r=1,$ subsequently proved by B\'ona \cite{MB98}. B\'ona also proved \cite{MB97} that the number of $r$-occurrences of the pattern 132 is P-recursive in the size. Equivalently, the ordinary generating function is D-finite. He then proved the stronger statement that the generating function is algebraic. 

Mansour and Vainshtein \cite{MV02} proved the corresponding result for $r=2$ and then gave conjectural results for $r \le 6,$ and conjectured the structure of the general formula.

For the increasing subsequence 123, Noonan \cite{JN96} proved the result for $r=1,$ and Noonan and Zeilberger \cite{NZ96} conjectured the result for $r=2.$ This was subsequently proved by Fulmek \cite{MF02}, who also gave conjectured results for $r=3$ and $r=4.$ These were subsequently proved by Callan \cite{DC02}. Nakamura and Zeilberger \cite{NZ13} developed a Maple package implementing a functional equation that readily generated terms for the cases $r \le 7,$ and gave the corresponding expressions for $\psi_r(n)$ for $r \le 7.$

In \cite{MB07} B\'ona proved that, for the monotone pattern $1234\cdots k,$ the distribution function of $r$-occurrences is asymptotically normal. Indeed, this is true for {\em any} classical pattern, a result first claimed to be proved by B\'ona in \cite{MB07a}\footnote{Mikl\'os B\'ona tellsus that this arXiv was never proceeded to publication, as he found an error in the proof.}. It was unequivocally proved by Janson et al. in \cite{JNZ13}, and can also be proved, perhaps even more easily, by the methods developed by Hofer in \cite{LH18}. 

 B\'ona, Sagan, and Vatter~\cite{bona:pattern-frequen:} looked at the sequences $\psi_r(n)$ for fixed $n$ and varying $r$ and for the pattern $132.$
They proved that there are both infinitely many values of $n$ for which this sequence has ``internal zeroes'' (so $\psi_r(n)=0$, but $\psi_q(n)>0$ and $\psi_s(n)>0$ for some $q<r<s$), and infinitely many values of $n$ for which it does not have internal zeroes.

In Brignall, Huczynska, and Vatter~\cite{brignall:decomposing-sim:}, the concept of simple permutations was introduced. It is known that a class with only finitely many ``simple permutations'' has an algebraic generating function. One of the simplest classes  with only finitely many simple permutations is $\Av(132)$. In Theorem 2.2 of  ~\cite{brignall:decomposing-sim:} they showed how this property extends to the class $\Av(132^{\le r})$ -- the set of permutations with at most~$r$ copies of $132$. (So permutations with precisely~$r$ copies of $132$ are just the difference of two algebraic generating functions).
 Theorem 2.3 of that paper shows that this algebraicity extends to counting the even permutations, involutions, even involutions, alternating permutations, and other permutations with at most~$r$ copies of $132$.

 As far as we are aware, the general situation for patterns of length 4 has not previously been studied, due in large part to the difficulty of generating coefficients of the underlying generating functions, though as we note below, there have been some series generated. We have developed an algorithm for this purpose, and find that of the $4!=24$ permutations of length-4, there are now seven effective Wilf classes. They are as follows:
\begin{eqnarray*}
I: & 1234, 4321\\
II:& 1243,2134,3421,4312,\\
III:& 1432, 2341, 3214, 4123\\
IV:& 2143, 3412,  \\
V:&1324, 4231\\
VI:&1342, 1423, 2314, 2431, 3124, 3241, 4132, 4213 \\
VII:& 2413, 3142
\end{eqnarray*}
By {\em effective Wilf class} in this context we mean that the number of $r$-occurrences of any pattern in the class in a permutation of length $n$ is the same, irrespective of the value of $r.$
For pattern-avoiding permutations, the first four classes above correspond to a single Wilf class, so these will all have coefficient growth $9^n.$ The fifth entry, with coefficient growth $\mu^n,$ where $\mu\approx 11.598 $ \cite{CGZ18} corresponds to another Wilf class, and the third Wilf class for PAPs comprises the patterns in class VI and VII above, with coefficient growth $8^n.$

It is worth remarking that, for contiguous PAPs, there are also 7 Wilf classes, and 5 of those 7 are the same as those given above. However the two problems differ in classes III and VI. As far as we can see, there is no reason that these two problems should be connected, so it's not surprising that there is a difference between the class distributions. It is perhaps surprising that 5 classes are the same.

\begin{figure}[ht!] 
\begin{minipage}[t]{0.45\textwidth} 
\centerline {\includegraphics[width=\textwidth]{ 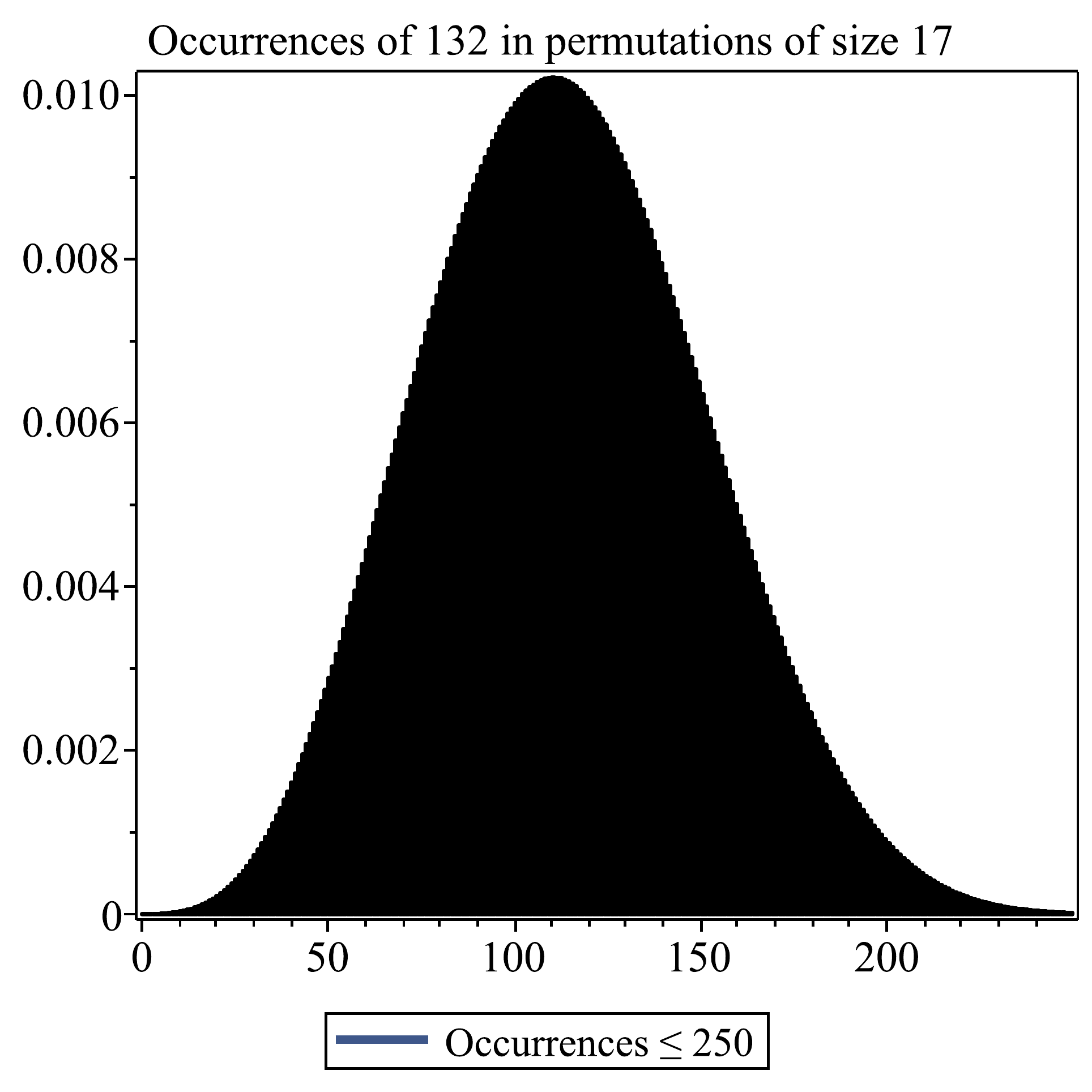}}
\caption{Histogram of  normalised 132 occurrences.} 
\label{fig:132}
\end{minipage}
\hspace{0.05\textwidth}
\begin{minipage}[t]{0.45\textwidth} 
\centerline{\includegraphics[width=\textwidth]{ 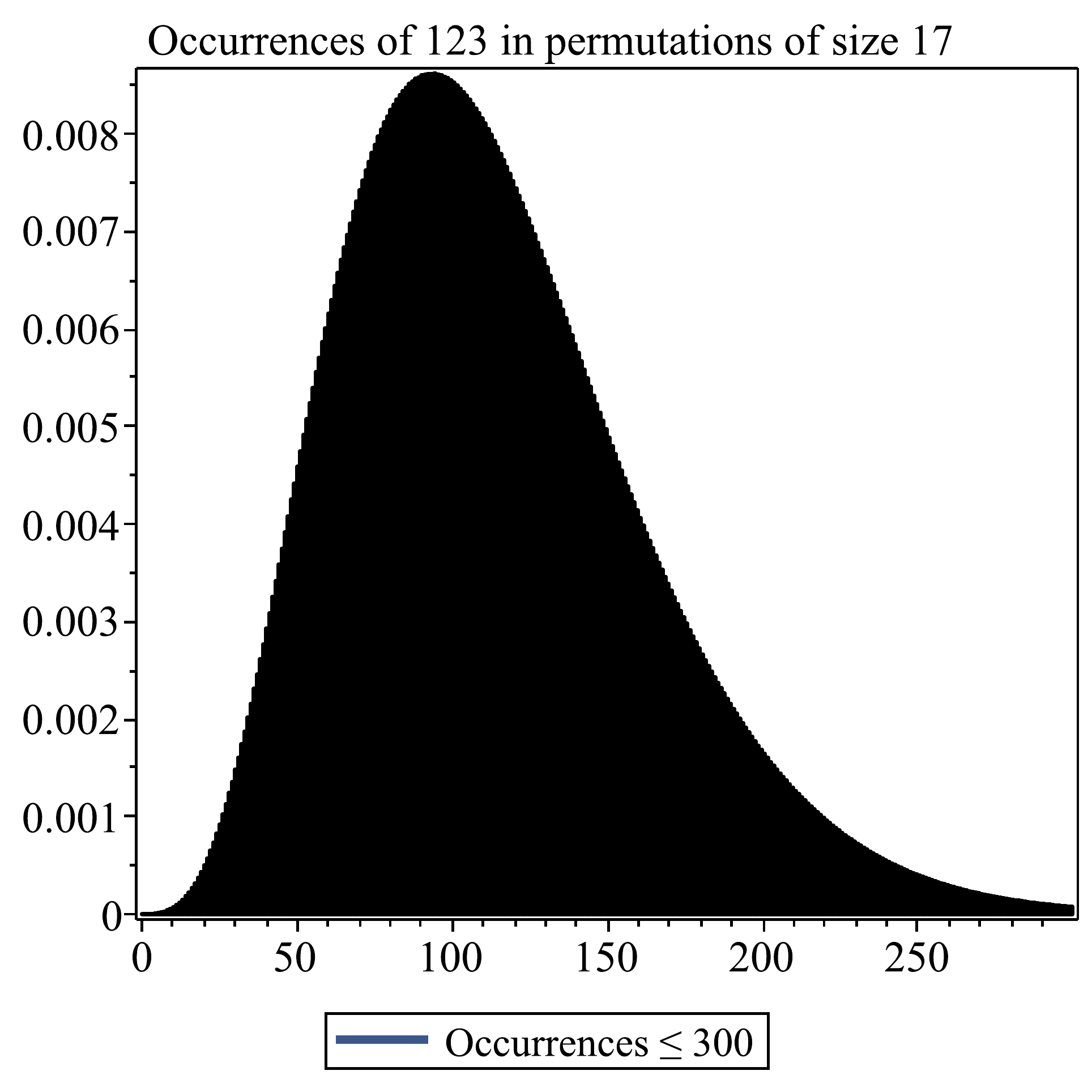}}
\caption{Histogram of normalised 123 occurrences.} 
\label{fig:123}
\end{minipage}
\end{figure}

\section{Generating the number of occurrences of a given pattern}

In \cite{M11} Minato came up with an innovative and general algorithm for counting pattern
avoiding permutations. It involved a straightforward algorithm operating upon sets of permutations
that generates all permutations containing the pattern, and then, crucially, an
efficient computational representation of said sets taking up space and time much smaller
than the number of elements in the set. This makes the algorithm much
more efficient in practice than other known general algorithms that examine
each pattern avoiding permutation individually. In \cite{I17}, Inoue
improved the set representation to make the algorithm even more efficient in practice.

We use the same algorithm operating upon sets, except apply it to multisets instead
of sets. That is, each element in the set has a multiplicity. The union of two multisets
in our context sums the multiplicities of each element. A consequence of the construction
algorithm is that the resulting set
contains each permutation containing the pattern, with a multiplicity equal to the
number of occurrences of the given pattern. It is then an efficient operation on the
multiset to obtain the number of elements with each multiplicity. This is the
desired enumeration.

We generalize the back end set representation in \cite{I17} to handle multisets. This
actually is a potentially generally useful multiset representation for combinatorics.
In practice, it ends up being reasonably efficient. As a more complex structure containing 
more information than the simple set it takes somewhat more time and memory than the simple set,
but is still vastly more efficient than techniques that consider elements individually.

The rest of this section goes into more detail, firstly on the set construction algorithm,
secondly on the method of representing a set of permutations as a binary decision tree,
and thirdly on generalising binary decision trees to contain multiplicities so as to be able to
represent multisets.

\subsection{Set construction algorithm}

Here we summarise the algorithm used in \cite{M11} and \cite{I17}, leaving out many
of the details and concentrating upon the set operations used.

For the following discussion, $n$ unless
otherwise specified is the number of elements in the desired permutation, and $k \leq n$ is the length
of the pattern. 

There is a three step process to generate the set of all pattern containing permutations.\footnote{
In \cite{M11}
and \cite{I17} this is followed by a trivial subtraction from the set of all permutations to get
pattern avoiding permutations; this subtraction is ignored in this paper as we are interested in
pattern containing permutations. 
} 

Step 1 generates all permutations $A$ whose first $k$ elements are in increasing order. There are
$\binom{n}{k}$ ways of choosing the first $k$ elements, and $(n-k)!$ ways of arranging the remaining
elements, for a total of $n!/k!$ elements in this set. We will not repeat the details of the construction
here, but mention that it generates each element exactly once and requires $O(n^2)$ set unions and 
set compositions with a ``unit'' permutation (a swap of two elements in \cite{M11} and a rotation in \cite{I17}).
If operating on a multiset, each element is created with multiplicity $1$.

As an example, if $n=4$ and $k=2$ for the pattern $21$ (only the length matters for this step), $A= \{1234, 1243, 1324, 1342, 1423, 1432, 2314, 2341, 2413, 2431,
3412, 3421 \}$ are the elements whose first $k=2$ elements are in order.

Step 2 makes a new set $B\times A$ that just rearranges the first $k$ elements of the above set $A$ to match
the pattern. That is, $B$ is a set containing a single element
which is a permutation consisting of the pattern followed
by the remaining $n-k$ elements of the identity permutation.

As an example, if $n=4$ and $k=2$ for the pattern $21$,
$B=\{ 2134 \} $.

The required set operation $\times$ is a cross product, taking two sets (or multisets) $A$ and $B$ and producing
a set $\{a\times b : a\in A, b\in B\}$ where $a\times b$ is the composition of the permutations $a$ and $b$.
If multisets are being used, then the cardinality (sum of the multiplicities of each element) of $A\times B$ is
the product of the cardinalities of $A$ and $B$.

In the running example, the only difference is that the
first two elements are exchanged to match the pattern $21$. Giving the elements in the same order as for $A$, 
$B \times A =  \{ 2134, 2143, 3124, 3142, 4123, 4132, 3214, 3241, 4213, 4231,
4312, 4321 \}$. 

$B$ has cardinality 1. 
The cardinality of $B\times A$ is the same as $A$: that is, $n!/k!$ elements each with multiplicity $1$.

Step 3 distributes the $k$ first elements of each permutation in $B\times A$ over the $\binom{n}{k}$ possible
positions in the permutation. This is done by constructing a set $C$ of permutations that take
the first $k$ elements to each of the $\binom{n}{k}$ possible
positions in the permutation. This construction is similar to the construction of $A$, requiring set
unions, and composition with unit permutations. The cardinality of $C$ is $\binom{n}{k}$ elements
each with multiplicity 1.

The set of pattern containing permutations is then $C\times B \times A$. The same pattern may
occur multiple times of course.

For $n=4$ and $k=2$, $C = \{ 1234, 1324, 1342, 3124, 3142, 3412 \}$.
Unfortunately $C\times B \times A$ has too many elements to helpfully list,
but as a few examples, the first permutation in $C$ is the
identity which when crossed with $B\times A$ changes nothing.
The second element in $C$, $1324$ (pattern positions $1$ and $3$), when operating on $2134$ (pattern values $2$ and $1$)
produces $2314$. Of course there is another instance of
that pattern in $C\times B \times A$: $3124$ (pattern positions $2$ and $3$) 
operating on $3124$ (pattern values $3$ and $1$) also
produces $2314$.

As a multiset, the number of ways the same pattern occurs in this construction for a permutation
$p\in C\times B\times A$ (the multiplicity of $p$) is the number of occurrences of the pattern
in $p$. This is exactly what is desired.\footnote{
As the multiset cardinality of $C\times B\times A$ is $\binom{n}{k} n!/k!$, this means that the
sum over all $n!$ permutations $p$ of the number of occurrences of the pattern is 
also $\binom{n}{k} n!/k!$, and so the mean number of occurrences of the pattern is
$\binom{n}{k}/k!$. This is a well known result, see for example Sec. 3.1 in \cite{JNZ13}.
}

\subsection{Set representation}

As mentioned before, the obvious simple computational set representation
explicitly listing all elements would result in dreadful
performance, as there are close to $n!$ elements in
$C\times B \times A$.

The big inspiration of \cite{M11} was to use binary 
decision diagrams (BDDs) to represent these sets\footnote{Actually
a slight variant of binary decision diagrams is used, zero 
suppressed binary decision diagrams (ZDDs). These have
similar computational properties, and their efficiency is within
a factor of their number of variables. In practice ZDDs are somewhat more efficient when most variables are false, such as the
present case. Hereafter we will describe from the viewpoint of BDDs 
which are conceptually slightly simpler; 
everything works with either BDDs or ZDDs. The
changes to ZDDs are straightforward, and should be
used in practice for this application.}.
A BDD represents a Boolean function of $N$ Boolean variables.
A Boolean function $f(v_0,...,v_{N-1})$ 
can be considered to represent a set $F$ with an element
$(v_0,...,v_{N-1})$ in $F$ iff $f(v_0,...,v_{N-1})$ is true.
For many combinatorial problems BDDs produce a very
compact representation of these sets \cite{Knuth09}.

A BDD consists of a list of $N$ Boolean variables ${v_0,...v_{N-1}}$,
a table consisting of some number of rows, and a starting index
pointing to one of these rows. Each row consists of three 
elements - a variable index, and two row indices (LO and HI). There are two
special row indices meaning true (often called 1 or $\top$) or 
false (often called 0 or $\bot$). One can store
multiple values in the table, in which case each starting index
is effectively a function. To evaluate a function implied by
a row index, first check to see if it is one of the two special
true or false references. If so, we are done. Otherwise read
the row corresponding to that index. Check the variable that
row references. If the variable is false, take the first
row index (LO) in the current row. Otherwise take the second row index (HI).
Replace the original row index with this new index, and 
repeat the above operation recursively. See figure \ref{fig:BDD_example}
for an example.

In order to promote uniqueness, the following rules are 
required:
\begin{itemize}
    \item If row $i$ contains a reference to a non-terminal row $j$,
    then the variable index in row $i$ must come before the variable 
    index in row $j$.
    \item No two rows in the table should be identical (if one were
    tempted to add a new row the same as an existing row, just use
    the existing row).
    \item Each row should matter. No row should have the same
    indices regardless of the variable state (except possibly the 
    special case of the two terminal rows, if they are explicitly stored).
\end{itemize}
This produces uniqueness (up to row renumbering) for any given
function, and, possibly unintuitively, produces a compact representation
for many combinatorics problems.

\begin{figure}
  \centering
  \begin{subfigure}[b]{0.3\textwidth}
    \centering
    \begin{tikzpicture}[grow = down]
      \node[](root){}
      child{
        node[draw,circle](A){0}
          child { node[draw,rectangle] { $\top$ } edge from parent[dashed] }
          child { 
            node[draw,circle]{1}
            child { node[draw,rectangle] { $\bot$ } edge from parent[dashed] }
            child { node[draw,rectangle] { $\top$ } }
          }
       };
    \end{tikzpicture}    
    \caption{Graphical representation. One starts from the
    top of the diagram. A circle represents a decision point.
    The number $i$ in the circle is the variable index. Take the
    LO choice (dashed) if $v_i$ is false and the  HI choice (solid) if $v_i$ is true.
    A square is a terminal : $\top$ means true and $\bot$ means false.}
  \end{subfigure}
  \hfill
  \begin{subfigure}[b]{0.3\textwidth}
    \centering
    \begin{tabular}{|cc|c|}
        \hline
        $v_0$ & $v_1$ & output   \\
        \hline
        false & false & true    \\
        false & true  & true    \\
        true  & false & false   \\
        true  & true  & true    \\
        \hline
    \end{tabular}
    \caption{Truth table. This shows the result of the
    BDD for all possible combinations of values of $v_0$ and $v_1$}
  \end{subfigure}
  \hfill
  \begin{subfigure}[b]{0.32\textwidth}
    \centering
    \begin{tabular}{|r|ccc|}
       \hline
       row & variable & LO & HI \\
       \hline
       3& 0 & 1 & 2  \\
       2& 1 & 0 & 1  \\
       1& \multicolumn{3}{c|}{special true $\top$} \\
       0& \multicolumn{3}{c|}{special false $\bot$} \\
       \hline
    \end{tabular}
    Start from row 3.
    \caption{BDD table. This is the computer representation of the
    BDD. Each row has an index (the row column). For each row one stores
    a variable index, and the LO and HI indices to other columns in this table.
    }
  \end{subfigure}

  \caption{Three different representations of a BDD of two variables.}
  \label{fig:BDD_example}
\end{figure}

An example BDD using notation similar to \cite{Knuth09} is given in 
figure \ref{fig:BDD_example}

ZDDs are very similar to BDDs, except instead of suppressing rows
whose two children are identical, one suppresses rows where a true
value of the variable would lead to the false special row index.
This ends up being more compact than a BDD when most variables are
false in most solutions.

The detailed mechanics of implementation of 
a BDD or ZDD are well described in \cite{Knuth09}.

Of course, we need a set of permutations, not a set of
choices of a set of binary variables. In \cite{M11}, Minato
defined a set of basis permutations - the $N=n(n-1)/2$
pairwise swaps of each element of the $n$ elements
of the permutation. Minato also defined a canonical order. Each
permutation then has a unique set of these $n(n-1)/2$
basis permutations that combine to produce it\footnote{
Of course there are $n!$ permutations and $2^{n(n-1)/2}$
choices of basis variables, so not all choices of the
Boolean variables are valid permutations. One has to be
careful not to accidentally construct one of these
invalid permutations in the set operations. For union and
intersection, there is no issue - if the inputs are valid,
the output will be. For the composition of two permutations,
the algorithms in \cite{M11} and \cite{I17} are designed
with this in mind.
}. A $\pi$DD is a  zero suppressed binary decision
diagram using this encoding to represent a set of
permutations.

The encoding is improved in \cite{I17} so as to have
a basis of $n(n-1)/2$ rotations of the elements between
each choice of two elements of the permutation. 
This variant of a  $\pi$DD is called a Rot-$\pi$DD.
This basis set
produces a much smaller representation of the set $A$ 
than a $\pi$DD which intuitively and empirically leads to a 
smaller representation of $C\times B\times A$, and a
more efficient algorithm.

The required set operations are set union, unit and
empty set creation, the cross product, and set
cardinality.

Unit and empty set creation are trivial and will
clearly not create an invalid permutation element
given a valid input. The set union and cardinality 
operations are identical to the standard BDD or ZDD 
operations \cite{Knuth09}, and will not create new (and 
therefore possibly invalid) permutation elements.

The cross product is not a standard BDD operation as 
it needs to take into account the composition of permutations
and the requirement for a canonical representation of
each permutation. In practice this is comparable to
the standard BDD or ZDD cross product, but is a little
more fiddly. Details are described in \cite{M11} and \cite{I17}.

There is another operation used in practice to create
the sets $A$ and $C$ which is really a special case of
the set cross product, where one of the sets is a single
element set containing a permutation which is 
one of the basis permutations. For efficiency
this is implemented as a special case. 
Details are described in \cite{M11} and \cite{I17}.

There are no useful theoretical bounds on the performance
of $\pi$DDs or Rot-$\pi$DDs for use in the above algorithm
to enumerate pattern avoiding permutations, 
but empirically they are vastly superior to algorithms
that consider each pattern avoiding permutation.

\subsection{Multiset representation}

Minato \cite{M11} and Inoue \cite{I17} use a ZDD as the data structure representing
a function of $N$ Boolean variables producing a Boolean answer,
where {\em true} means set membership. We use a similar data structure
to represent a function of $N$ Boolean variables producing a
non-negative integer which represents the multiplicity of the
element in the set.

The main change to the BDD (or ZDD) data structure needed 
is to replace the concept of a row reference $r$ by a tuple $(r,m)$
of a row reference $r$ and a natural number multiplier $m$. A row
now contains a variable index and two of these tuples, instead of 
just a variable index and two row indices.\footnote{
There are other ways BDDs could be generalized to multisets. For instance
one could have just one multiplier associated with each row, rather
than one multiplier for each of the two children. This seems sensible
as it would store less per row. However, it would mean having a 
larger number of rows, one for each different multiple ever associated
with such a row. It also makes it somewhat messy to store a multiple
of the special terminal row representing true. That is, it would
be messy to represent the set containing all elements, each with
multiplicity two}

A multiset version of a BDD we call an MBDD; a multiset version
of a ZDD we call an MZDD.

Evaluation of the function is the same as for BDDs, although
all the multipliers $m$ on the taken path need to be multiplied
together to get the final multiplicity. The special terminal {\em false}
is considered to have a multiplicity of $0$, the special terminal true
has a multiplicity of $1$ (times, of course, any multipliers {\em en route}).

In order to maintain uniqueness of representation (and thus
compactness and thus efficiency), one has to add two new rules:
\begin{itemize}
    \item in any row, the greatest common divisor of the two multipliers in that row must be one. If one is tempted to want a row with a larger common
    divisor, factor out the gcd, and include it in the multiplier for the
    reference to the row.
    \item  The multiplier associated with a row reference
    for the special terminal {\em false} should always be zero.\footnote{
    This means one could use a single terminal instead of two terminals,
    resolved using the multiplicity. The efficiency gain from this is
    tiny, and it complicates the analogy to normal BDDs.
    }  
\end{itemize}

\begin{figure}
  \centering
  \begin{subfigure}[b]{0.28\textwidth}
    \centering
    \begin{tikzpicture}[grow = down]
      \node[](root){}
      child{
        node[draw,circle](A){0}
          child { node[draw,rectangle] { $\top$ } edge from parent[dashed] node[left]{7}}
          child { 
            node[draw,circle]{1} 
            child { node[draw,rectangle] { $\bot$ } edge from parent[dashed] node[left]{0} }
            child { node[draw,rectangle] { $\top$ } edge from parent node[right]{1}}
            edge from parent node[right]{13}
          }
          edge from parent node[right]{2}
       };
    \end{tikzpicture}    
    \caption{Graphical representation. Interpretation is the
    same as figure \ref{fig:BDD_example}.
    The result is the product of the numbers labeling the edges taken.}
  \end{subfigure}
  \hfill
  \begin{subfigure}[b]{0.37\textwidth}
    \centering
    \begin{tabular}{|cc|c|}
        \hline
        $v_0$ & $v_1$ & output   \\
          &   & multiplicity   \\
        \hline
        false & false & $2\times 7=14$    \\
        false & true  & $2\times 7=14$    \\
        true  & false & $2\times 13\times 0=0$   \\
        true  & true  & $2\times 13\times 1=26$    \\
        \hline
    \end{tabular}
    \caption{Multiplicity table. Like the truth table in 
    figure \ref{fig:BDD_example} except the output is a multiplicity.}
  \end{subfigure}
  \hfill
  \begin{subfigure}[b]{0.31\textwidth}
    \centering
    \begin{tabular}{|r|ccccc|}
       \hline
       row & v & \multicolumn{2}{c}{LO} & \multicolumn{2}{c|}{HI} \\
           &          & $r$ & $m$ & $r$ & $m$ \\
       \hline
       3& 0 & 1 & 7 & 2 & 13  \\
       2& 1 & 0 & 0 & 1 & 1  \\
       1& \multicolumn{5}{c|}{special true $\top$} \\
       0& \multicolumn{5}{c|}{special false $\bot$} \\
       \hline
    \end{tabular}
    Start from row 3, multiple 2.
    \caption{BDD table. The LO and HI fields now have
    a multiplicity $m$ in addition to the row index $r$.
    }
  \end{subfigure}

  \caption{Three different representations of a MBDD of two variables, 
  similar to figure \ref{fig:BDD_example}.}
  \label{fig:MBDD_example}
\end{figure}
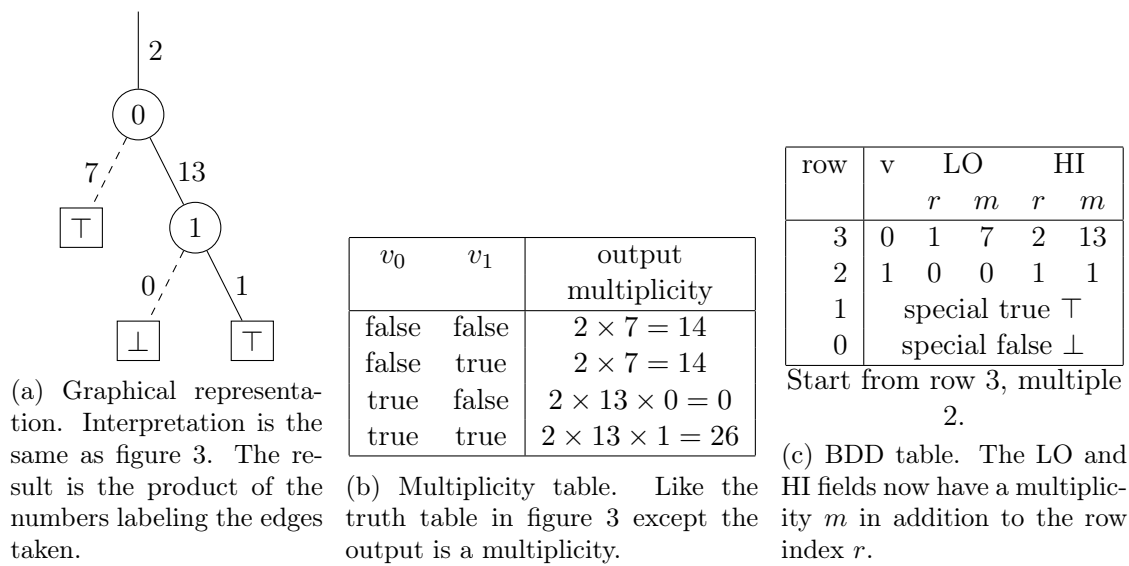

Figure \ref{fig:MBDD_example} is an example of an MBDD analogous to 
figure \ref{fig:BDD_example}.

In practice, this representation works well, and it is very straightforward 
to generalize the BDD or ZDD set union and cross product algorithms 
(and many other BDD algorithms) 
to include the multiples\footnote{Explicit details, if required, are
given in the software available below. The only changes are to track
multiples - adding them when doing unions, and multiplying them for 
intersections, and resolving gcd constraints.}.

\subsubsection{Cardinality algorithm of a MBDD}
The cardinality algorithm on BDDs or ZDDs (Algorithm C in \cite{Knuth09} section 7.1.4) can also be generalized in
a straightforward manner, multiplying each $c_i$ by the associated multiple, 
although now one gets the sum of the
multiplicities of all elements rather than the number of elements.
This is not however the operation we will wish to do - we want to
get the number of elements with each multiplicity in the data structure.

As in the standard cardinality algorithm, we assume that the table
is sorted topologically. That is, the LO and HI fields of each row
point to lower rows. This property drops naturally out of the construction
algorithm and doesn't have any computational cost.

What we want is the number of elements in the set (otherwise known
as combinations of the input variables) that have a given multiplicity.
Represent this as an array of length equal to the highest multiplicity
plus one\footnote{The plus one is to account for the multiplicity
zero term. The algorithm is described including this term, although
in practice one can leave it out as it is redundant - the sum of the
array must equal two to the power of the number of variables considered. Leaving it
out produces a minor reduction in memory use due to the arrays being smaller,
and, as it is often the largest value by a long way, 
means the integers in the arrays will be smaller allowing fewer bits to be
used in their representation, saving a significant amount of memory.}.
We will call this array a generating function.
Define $G(r)$ to be the generating function for row $r$ considering
all variables from the variable in row $r$ and below.

Like the standard cardinality algorithm, we will calculate a table
for each row of $G(r)$ starting from the special terminal rows and 
working upwards.

Define $H(r,m,i)$ as the generating function represented by row $r$
with extra multiple $m$ considering variables $v_i$ and after.
Then $H(r,m,0)$ is the desired generating function for the starting
row $r$ and multiple $m$. In particular, for the example in figure~\ref{fig:MBDD_example},
$H(3,2,0)$ is an array with 1 at index (multiplicity) 0, 2 at index 14,
and 1 at index 26, and zero elsewhere.

$H(r,m,i)$ can be easily computed from $G(r)$. Let $G(r)_j$ be the 
value of $G(r)$ with index $j$, and use the same subscript notation
for $H(r,m,i)$. Let $V_r$ be the variable index in row $r$, or
the total number of variables $V$ if $r$ is one of the special terminal rows.
Then
$$H(r,m,i)_{mj}=2^{V_r-i} G(r)_j$$
with all other values zero. Note that the length of $H(r,m,i)$ will
be $m$ times the length of $G(r)$, not counting the 0 index element
which generally is not needed in practice. $H(r,m,i)$ is not meaningful
if $i>V_r$.

Note that throughout this we maintain the invariants that the
sum of elements of $G(r)$ is $2^{V-V_r}$ and the sum of elements
in $H(r,m,i)$ is $2^{V-i}$.

For the terminals, $G(0=\bot)=\{1\}$ represents one instance 
with multiplicity zero. $G(1=\top)=\{0,1\}$ represents zero instances
with multiplicity zero and one instance with multiplicity one.

For each non-terminal,
$$G(r) = H(LO(r)_r,LO(r)_m,V_{r}+1) + H(HI(r)_r,HI(r)_m,V_{r}+1)$$
where $LO(r)_r$ is the row reference in the LO field of row $r$,
and $LO(r)_m$ is the multiplicity in the LO field of row $r$,
and similar definitions for $HI(r)$ being the HI field of row $r$,
and addition of two generating functions is just elementwise addition
of corresponding indices.

The topological sorting means that by computing $G(r)$ in ascending
order, the right hand side only references already computed values of
$G$.

This is very similar to the standard BDD cardinality algorithm;
the same changes can be applied to the standard ZDD cardinality
algorithm to compute $G(r)$ for a MZDD. The only difference
is that the adjustment for variable gaps is different as non-consecutive
variables means a large number of zero multiplicity solutions:
\begin{equation}
  H(r,m,i)_{mj}=G(r)_j + \begin{cases}
    2^{V-i}-2^{V-V_r}, & \text{if $j=0$}.\\
    0, & \text{otherwise}.
  \end{cases}
\end{equation}

\subsection{Software availability}

A general library implementing BDDs, ZDDs, $\pi$DDs, 
and Rot-$\pi$DDs with either sets or multisets 
in the Rust language is
available at \url{ https://github.com/AndrewConway/xdd}. 
Inside this repository is an example program, {\tt examples/pap.rs}
that computes the number of permutations of length up to $n$,
split up by the number of times a given
pattern is contained. This is the program used for the
enumerations in this paper.

The library is also available on crates.io as {\tt xdd}.

\section{Counting occurrences of 132, 231, 312, 213}
Let $\psi_r(n)$ denote the number of permutations of length $n$ containing exactly $r$ occurrences of the nominated pattern.
Let $\Psi_r(x)$ be the ordinary generating function for $\psi_r(n).$\\
The generating function was shown by B\'ona \cite{MB97} to behave as
\BE \label{eqn:132}
\Psi_r(x)=\frac{1}{2} \left ( P_1(x)+P_2(x)(1-4x)^{-r+1/2}\right ),
\EE
where $P_1(x)$ and $P_2(x)$ are polynomials. Results are given for $r \le 5$ in \cite{MV02} by Mansour and Vainshtein. They find, where we have explicitly only shown results for $r < 4,$

$$\psi_0(n) = \frac{(2n)!}{n!(n+1)!}$$
$$\psi_1(n) = \frac{(2n-3)!}{n!(n-3)!}$$
$$\psi_2(n) = \frac{(n^3+17n^2-80n+80)(2n-6)!}{n!2!(n-4)!}$$
$$\psi_3(n) = \frac{(n^6+51n^5-407n^4-99n^3+7750n^2-22416n+20160)(2n-9)!}{n!3!(n-5)!},$$ and more generally,
$$\psi_r(n)=\frac{Q_r(n)(2n-3r)!}{n!r!(n-r-2)!} \sim \frac{4^{n-3r/2} \cdot n^{r-3/2}}{r!\sqrt{\pi}}$$
where $Q_r(n)$ is a polynomial of degree $3(r-1)$  whose leading-order term is precisely $n^{3(r-1)}.$

From the results in \cite{MV02} it appears that the polynomials in eqn.(\ref{eqn:132}), $P_1(x)$ and $P_2(x),$ are 
polynomials with integer coefficients, of degree $r$ and $2r+1$ respectively.
 
In fig. \ref{fig:132} we show the histogram of $\psi_r(17)/17!$ for the pattern $132.$ Note that $\sum_r \psi_r(17)/17!=1.$ There are in fact further entries out to $r=342,$ but we have truncated the histogram at $r=250$ as subsequent contributions are visually indistinguishable from 0. In fig. \ref{fig:123} we show the corresponding histogram for the pattern $123.$  As we mentioned previously, these distributions have been proved to be asymptotically normal \cite{JNZ13}, and we see hints of this in fig. \ref{fig:132} (if we ignore the long tail), but there is still considerable skewness visible in fig. \ref{fig:123}, so we conclude that permutations of size 17 are far from the asymptotic regime.

\section{Counting occurrences of 123 and 321}
Let $C(x) =\frac{1}{2}(1-\sqrt{1-4x})$ be the generating function for Catalan numbers. Then it is well known that
$$\psi_0(n) = \frac{(2n)!}{(n+1)!n!}\sim \frac{ 4^n}{n\sqrt{n\pi}}$$
$$\Psi_0(x)=\frac{1}{2x}(1-\sqrt{1-4x}).$$
With 1 occurrence of the pattern, Noonan proved in \cite{JN96}, that
$$\psi_1(n) =  \frac{3}{n} \binom{2n}{n-3}=\frac{6(2n-1)!}{(n+3)!(n-3)!}=\frac{6(n^2-3n+2)(2n-1)!}{(n+3)!(n-1)!}\sim \frac{3\cdot 4^n}{n\sqrt{n\pi}},$$
From this we derive $$\Psi_1(x)=\frac{1}{2x^3}\left ((1-6x+9x^2-2x^3)-(1-4x+3x^2)\sqrt{1-4x}\right ).$$
For 2 occurrences, Noonan and Zeilberger \cite{NZ96} conjectured, and Fulmek \cite{MF02} has proved, that
\begin{align*}
    &\psi_2(n) = \frac{(59n^2+117n+100)(2n-2)!}{(n+5)!(n-4)!}\\
    &= \frac{(n^2-5n+6)(59n^2+117n+100)(2n-2)!}{(n+5)!(n-2)!}\sim \frac{59\cdot 4^n}{4n\sqrt{n\pi}}.
\end{align*}
From this we obtain
\begin{align*}
    &\Psi_2(x)=\frac{1}{2x^5}(1-8x+20x^2-17x^3+7x^4-5x^5)\\
    &-\frac{1}{2x^5}(1-6x+10x^2-5x^3+3x^4-x^5)\sqrt{1-4x}.
\end{align*}

From data generated by our program we find
\begin{align*}
\psi_3(n) = &\frac{4n(113n^3+506n^2+937n+1804)(2n-3)!}{(n+7)!(n-5)!}\\
 =&\frac{4n(n^2-7n+12)(113n^3+506n^2+937n+1804)(2n-3)!}{(n+7)!(n-3)!} \sim \frac{113\cdot 4^n}{2n\sqrt{n\pi}}
\end{align*}
from which we find
$$\Psi_3(x)=\frac{1}{2x^7}\left (P_1(x)-P_2(x)\sqrt{1-4x}\right ),$$
where $$P_1(x)=1-10x+33x^2-32x^3-31x^4+70x^5-35x^6+2x^8,$$ and $$P_2(x)=1-8x+19x^2-6x^3-27x^4+28x^5-7x^6-2x^7.$$

For 4 occurrences we find
$$\psi_4(n) = \frac{P_8(n)(2n-4)!}{(n+9)!(n-4)!} \sim \frac{3561\cdot 4^n}{16n\sqrt{n\pi}},$$
where
\begin{align*}
P_8(n)=&3(1187n^8+1042n^7-15602n^6+4128n^5-219697n^4+876878n^3\\
&+1840192n^2+8761152n-13063680).
\end{align*}
For the generating function, we find
$$\Psi_4(x)=\frac{1}{2x^9}\left (P_1(x)-P_2(x)\sqrt{1-4x}\right ),$$
where 
$$P_1(x)=1-12x+50x^2-65x^3-107x^4+437x^5-588x^6+492x^7-314x^8+108x^9-3x^{10},$$ and $$P_2(x)=1-10x+32x^2-17x^3-107x^4+245x^5-256x^6+192x^7-102x^8+18x^9+x^{10}.$$
We subsequently found that the formulae for $\Psi_3(n)$ and $\Psi_4(n)$ (but not for  $\psi_3(n)$ and $\psi_4(n)$) had previously been conjectured by Fulmek \cite{MF02}. Similarly, the results we give below for $\psi_r(n)$ have been obtained previously by Nakamura and Zeilberger \cite{NZ13}, but the expressions for the generating functions $\Psi_r(x)$ for $r>4$  are believed to be new.

For 5 occurrences of the pattern, 
$$\psi_5(n) = \frac{P_{10}(n)(2n-5)!}{(n+11)!(n-5)!} \sim \frac{13123\cdot 4^n}{16n\sqrt{n\pi}},$$
where
\begin{align*}
P_{10}(n)=&2(-16445721600+1805846400n+1400051832n^2+803816948n^3\\
&+190602806n^4+35652267n^5-4824225n^6+11262n^7-57936n^8\\
&+68323n^9+13123n^{10}),
\end{align*}
and the generating function is
$$\Psi_5(x)=\frac{1}{2x^{11}}\left (P_1(x)-P_2(x)\sqrt{1-4x}\right ),$$
where 
\begin{align*}
P_1(x)&=1-14x+71x^2-126x^3-176x^4+1160x^5-2167x^6+2282x^7-1976x^8\\
&+1902x^9-1608x^{10}+824x^{11}-153x^{12}-2x^{13},\\
P_2(x)&=1-12x+49x^2-48x^3-212x^4+744x^5-1057x^6+956x^7-860x^8+838x^9\\
&-600x^{10}+212x^{11}-13x^{12}.
\end{align*}

For 6 occurrences of the pattern, 
$$\psi_6(n) = \frac{P_{12}(n)(2n-6)!}{(n+13)!(n-6)!} \sim \frac{193311\cdot 4^n}{64n\sqrt{n\pi}},$$
where
\begin{align*}
P_{12}(n)=&3(-9392423040000-2303171015040n+427865125056n^2+341897564488n^3\\
&+104336894932n^4+20001951630n^5+4778554443n^6+963184194n^7\\
&+135476931n^8+31717010n^9+4345001n^{10}+783318n^{11}+64437n^{12}),
\end{align*}
and the generating function is
$$\Psi_6(x)=\frac{1}{2x^{13}}\left (P_1(x)+P_2(x)\sqrt{1-4x}\right ),$$
where 
\begin{align*}
P_1(x)&=1-16x+96x^2-223x^3-192x^4+2295x^5-5493x^6+6299x^7-3491x^8\\
&+1098x^9-2070x^{10}+4777x^{11}-6187x^{12}+4525x^{13}-1486x^{14}+93x^{15},\\
P_2(x)&=-1+14x-70x^2+107x^3+312x^4-1625x^5+2903x^6-2473x^7+925x^8\\
&-436x^9+1398x^{10}-2581x^{11}+2687x^{12}-1439x^{13}+260x^{14}+x^{15}.
\end{align*}

Nakamura and Zeilberger \cite{NZ13} also give $\psi_7(n).$

The general situation seems to be
$$\psi_r(n) = \frac{Q_{2r}(n)(2n-r)!}{(n+2r+1)!(n-r)!} \sim \frac{C_r\cdot 4^n}{n\sqrt{n\pi}}$$
where $Q_{2r}(n)$ is a polynomial with integer coefficients of degree $2r.$
The amplitude coefficient $C_r$ appears to increase exponentially, growing seemingly like $\lambda^r,$ where the growth constant
$\lambda \approx 2.67.$ That is to say, the asymptotics remain unchanged, and only the amplitude, or premultiplying constant changes with $r.$ However this exponential increase in the amplitude cannot continue indefinitely, as the histogram plotting $\psi_r(n)$ against $r$ is asymptotically normally distributed. So given that, asymptotically, only the amplitude changes as $r$ changes, it must first increase and then decrease, reflecting the heights of the various histogram entries.

The generating function is conjectured to behave as
$$\Psi_r(x)=\frac{1}{2x^{2r+1}}\left (P_1(x)-P_2(x)\sqrt{1-4x}\right ),$$
where $P_1(x)$ and $P_2(x)$ are polynomials with integer coefficients whose degree depends on the parity of $r.$
If $r$ is even, both polynomials are of degree $5r/2.$ If $r$ is odd, $P_1(x)$ is of degree $(5r+1)/2,$ and $P_2(x)$ is of degree $(5r-1)/2.$

The conjectured form of the generating function has also been given previously by Fulmek \cite{MF02}, though without comment on the degree of the polynomials.

In fig. \ref{fig:123} we show the histogram of $\psi_r(17)/17!.$ Note that $\sum_r \psi_r(17)/17!=1.$ There are in fact further entries out to $r=680,$ but we have truncated the histogram at $r=300$ as subsequent contributions are visually indistinguishable from 0.

\section{Length 4 patterns}
In this section we investigate some properties of length-4 sequences. We have been able to generate data for permutations up to size 14, for all values of $r,$ but to go further requires greater computing resources than we have. The principal limitation is memory. We had 2TB at our disposal, but even that is insufficient to go beyond $n=14.$ As a consequence, we have been unable to conjecture any exact results for $\psi_r(n)$ for $r >0$ for any pattern, though we have been able to conjecture quite a lot about the asymptotics.

Also, if we restrict ourselves to $r=1$ and $r=2$ we can go further. It is also the case that specialising to a particular Wilf class sometimes allows a specialised algorithm to be developed that is more efficient than our general-purpose algorithm. In particular, we point out that for Wilf class I, Nakamura and Zeilberger \cite{NZ13} developed a functional equation approach that allowed them to obtain 70 terms for this class for $r=1,$ and these can be found as sequence A217057 in the OEIS \cite{OEIS}. Furthermore, Kauers pushed this recurrence to generate 200 terms, and these can be found at \cite{MK}.
For $r=2,$ these same authors provide 25 terms as entry A224249 in the OEIS \cite{OEIS}.

For Wilf class II, Nakamura \cite{BN13} gives 25 terms for $r=1,$ and these can be found as sequence A224179 in \cite{OEIS}, and 24 terms for $r=2$ as sequence A224249 in \cite{OEIS}.

For Wilf class III, Nakamura \cite{BN13} gives 18 terms for $r=1,$ and these can be found as sequence A224182 in \cite{OEIS}. For $r=2$ there appears to be no pre-existing data. We find 
\begin{align*}
&0,0,0,0,5,68,626,5038,38541,289785,2172387,\\
&16339840,123650958,942437531,7236542705.
\end{align*}
for the first 15 coefficients. That is to say, the $i$th coefficient is the number of permutations of length $i$ containing precisely two occurrences of the given pattern. Our algorithm would require more than 2TB of memory to go further.

For Wilf class IV there is no pre-existing data for $r>0,$ and the relevant sequences are not to be found in the OEIS. We find

For $r=1,$
\begin{align*}
  & 0,0,0,1,11,88,642,4567,32443,232189,1679295,12282794,90834993,\\
  &678779256,5121534664,38988595387,299244027539,2314045427659.\\
  \end{align*}
  
  And for $r=2,$
  \begin{align*}
  &0,0,0,0,4,53,495,4099,32345,250371,1926145,14820037,\\
  &114394941,887176357,6917420887, 54237535517.
  \notag
\end{align*}

For Wilf class V, Johansson and Nakamura \cite{JN13} developed the functional equations for the case of general $r,$ and in OEIS sequence A224182 they give the first 17 terms for the case $r=1,$ but give no values for $r=2,$ though they present the machinery for doing so.
We give the first 15 coefficients for $r=2$ in this case:
\begin{align*}
&0, 0, 0, 0, 6, 74, 645, 5023, 37549, 277089, \\
&2043416, 15146147, 113147663, 852978562, 6492322934.\\
\end{align*}

\begin{figure}[ht!] 
\begin{minipage}[t]{0.45\textwidth} 
\centerline{\includegraphics[width=\textwidth]{ 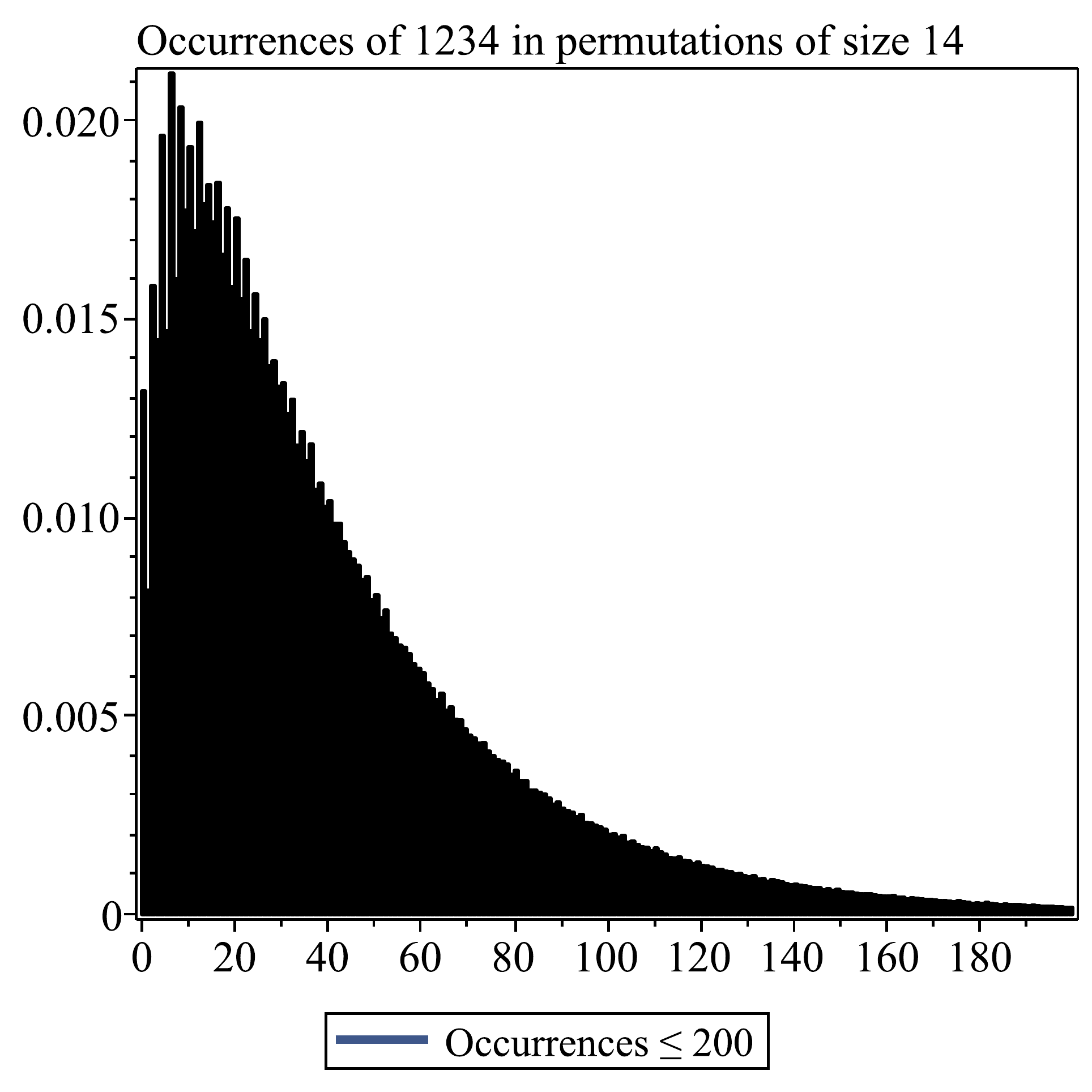}}
\caption{Histogram of 1234 occurrences.} 
\label{fig:1234}
\end{minipage}
\hspace{0.05\textwidth}
\begin{minipage}[t]{0.45\textwidth} 
\centerline {\includegraphics[width=\textwidth]{ 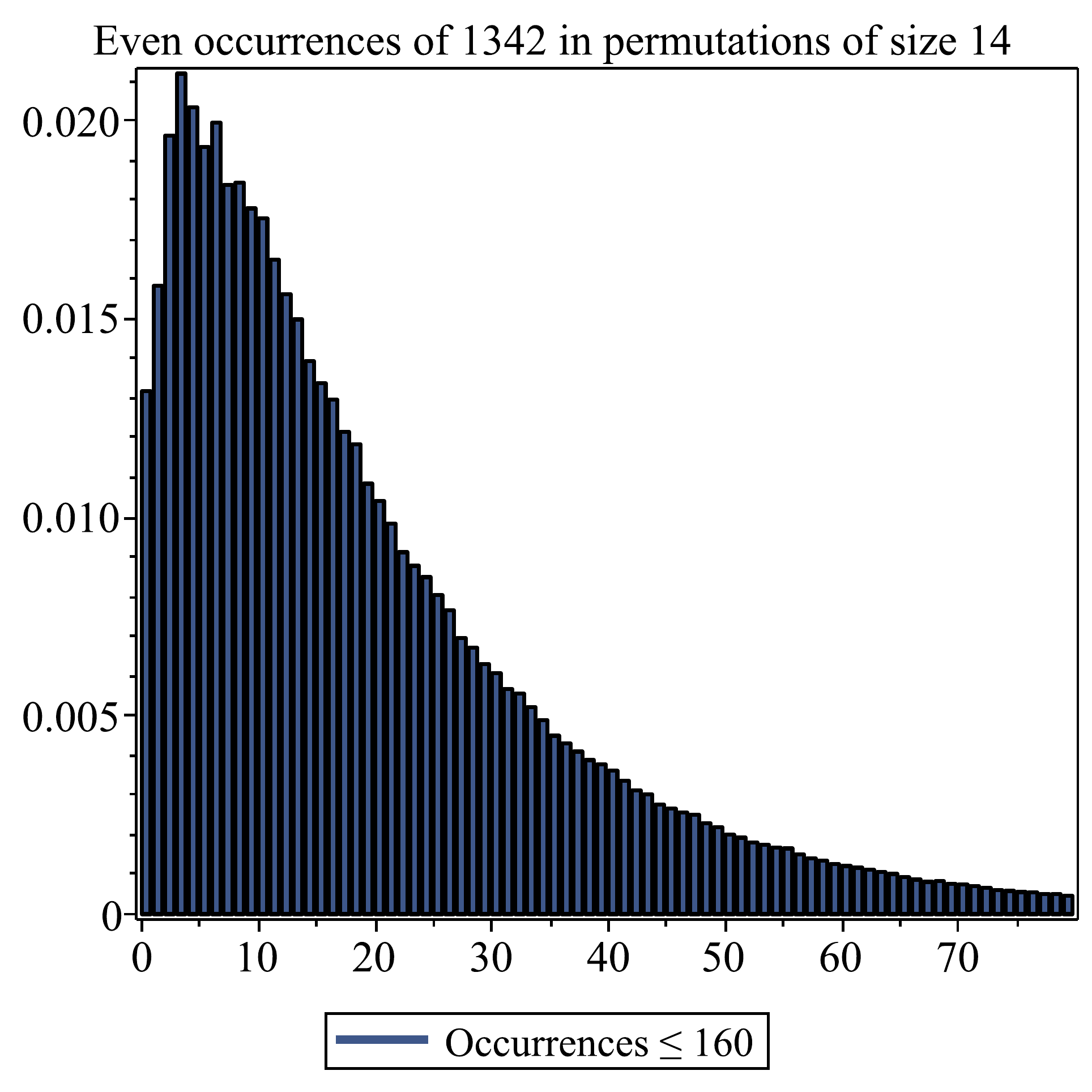}}
\caption{Histogram of 1234 even occurrences.} 
\label{fig:1234e}
\end{minipage}
\end{figure}

\begin{figure}[ht!] 
\begin{minipage}[t]{0.45\textwidth} 
\centerline{\includegraphics[width=\textwidth]{ 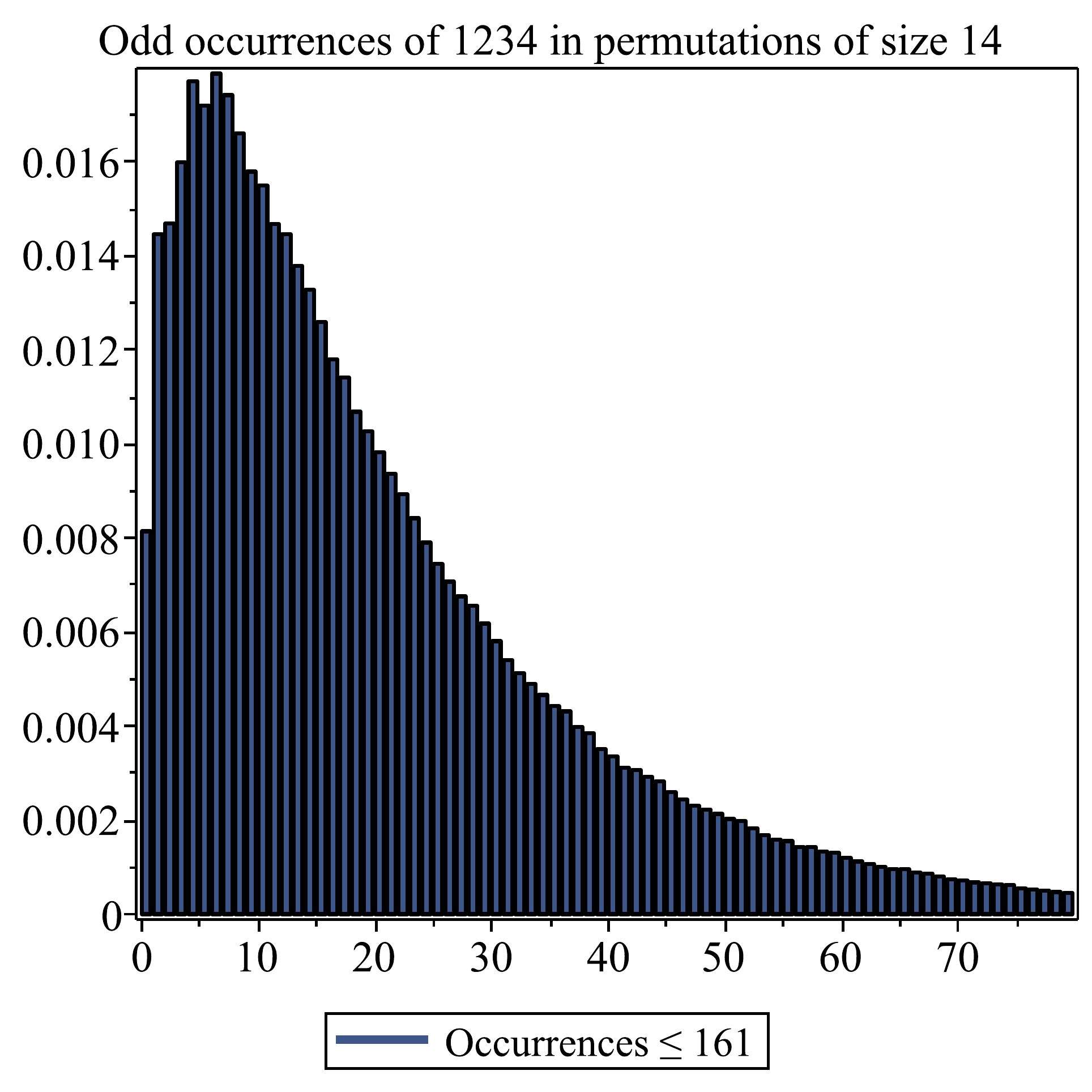}}
\caption{Histogram of 1234 odd occurrences.} 
\label{fig:1234o}
\end{minipage}
\hspace{0.05\textwidth}
\begin{minipage}[t]{0.45\textwidth} 
\centerline {\includegraphics[width=\textwidth]{ 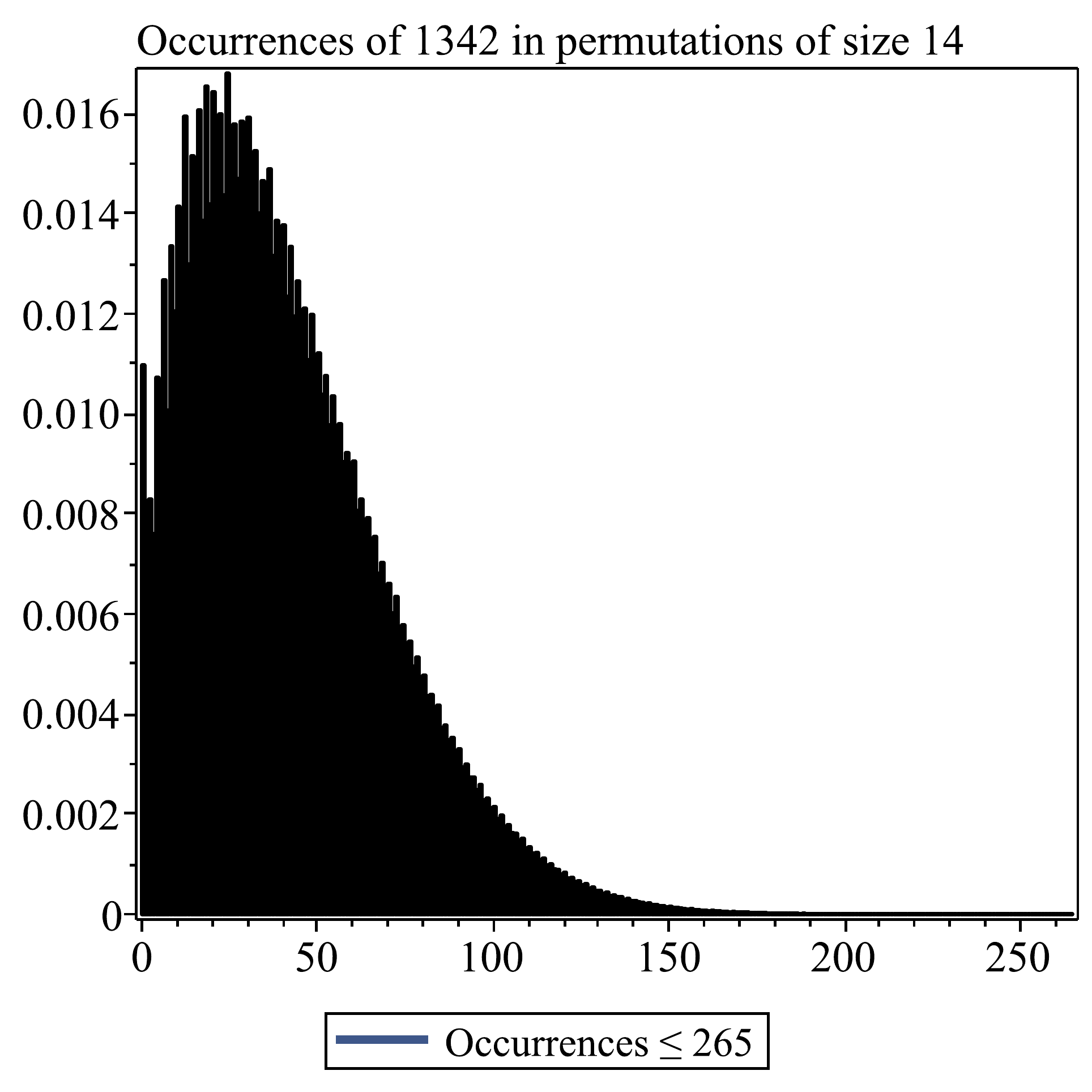}}
\caption{Histogram of 1342 occurrences.} 
\label{fig:1342}
\end{minipage}
\end{figure}

\begin{figure}[ht!] 
\begin{minipage}[t]{0.45\textwidth} 
\centerline{\includegraphics[width=\textwidth]{ 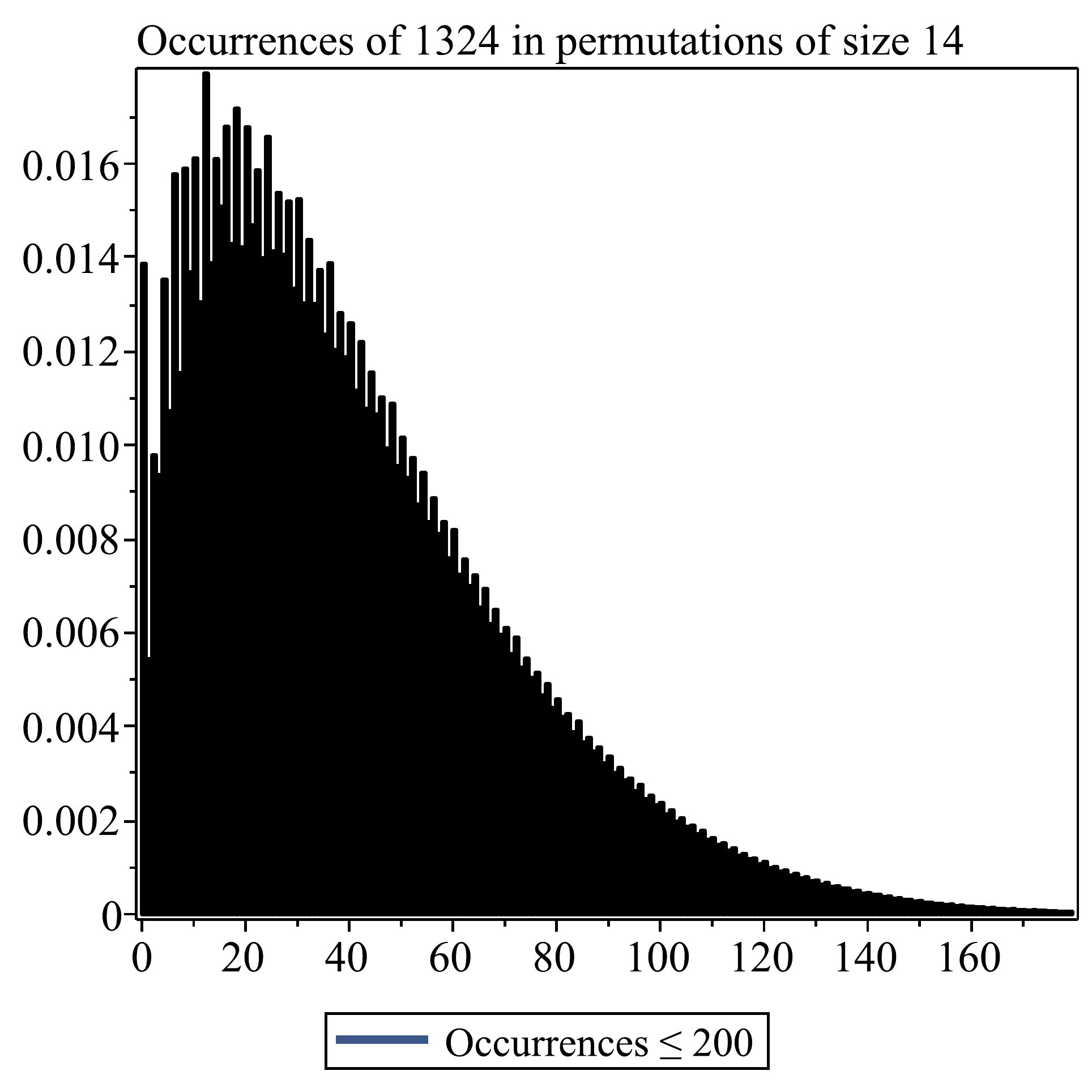}}
\caption{Histogram of 1324 odd occurrences.} 
\label{fig:1324}
\end{minipage}
\hspace{0.05\textwidth}
\begin{minipage}[t]{0.45\textwidth} 
\centerline {\includegraphics[width=\textwidth]{ 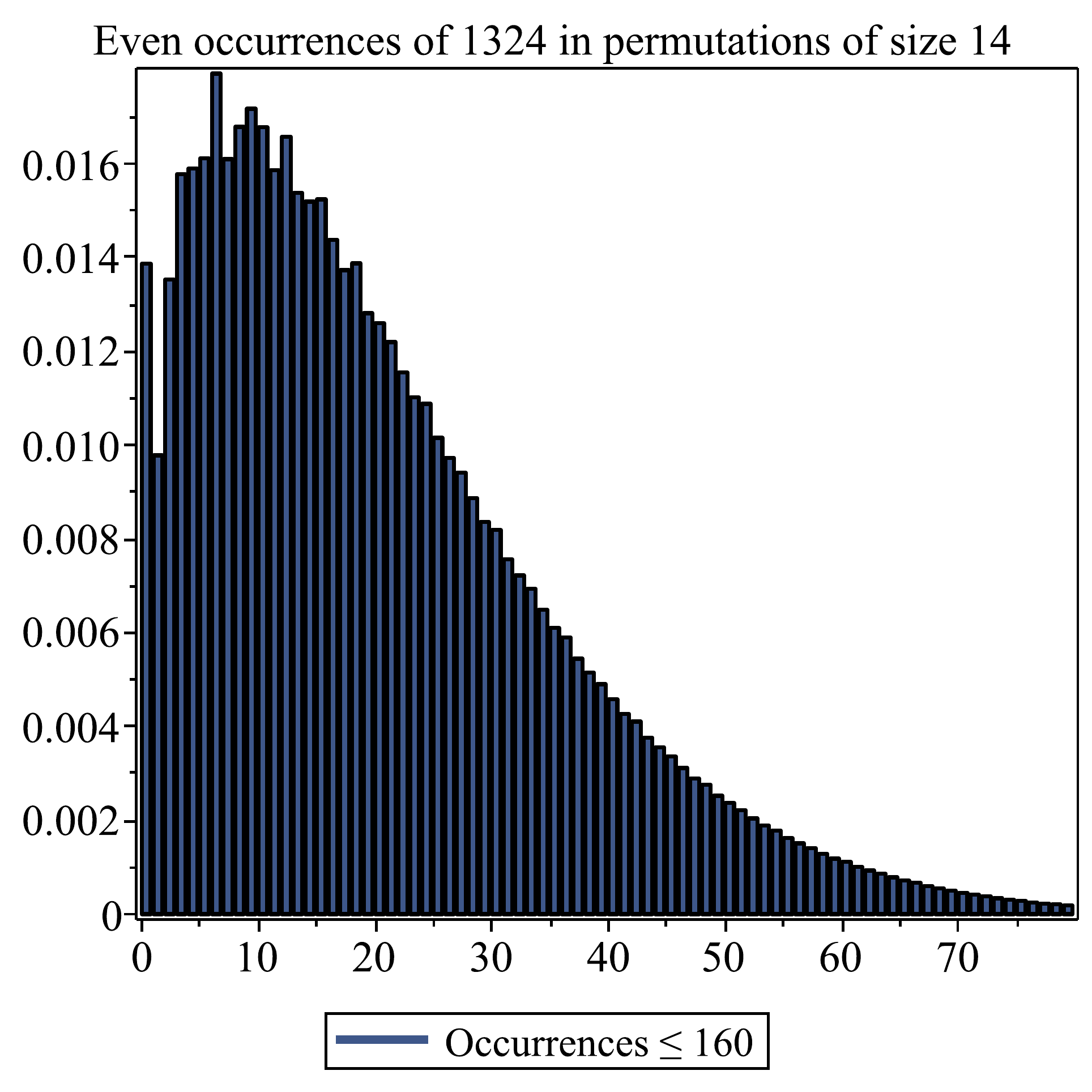}}
\caption{Histogram of 1324 even occurrences.} 
\label{fig:1324e}
\end{minipage}
\end{figure}

\begin{figure}[ht!] 
\begin{minipage}[t]{0.45\textwidth} 
\centerline{\includegraphics[width=\textwidth]{ 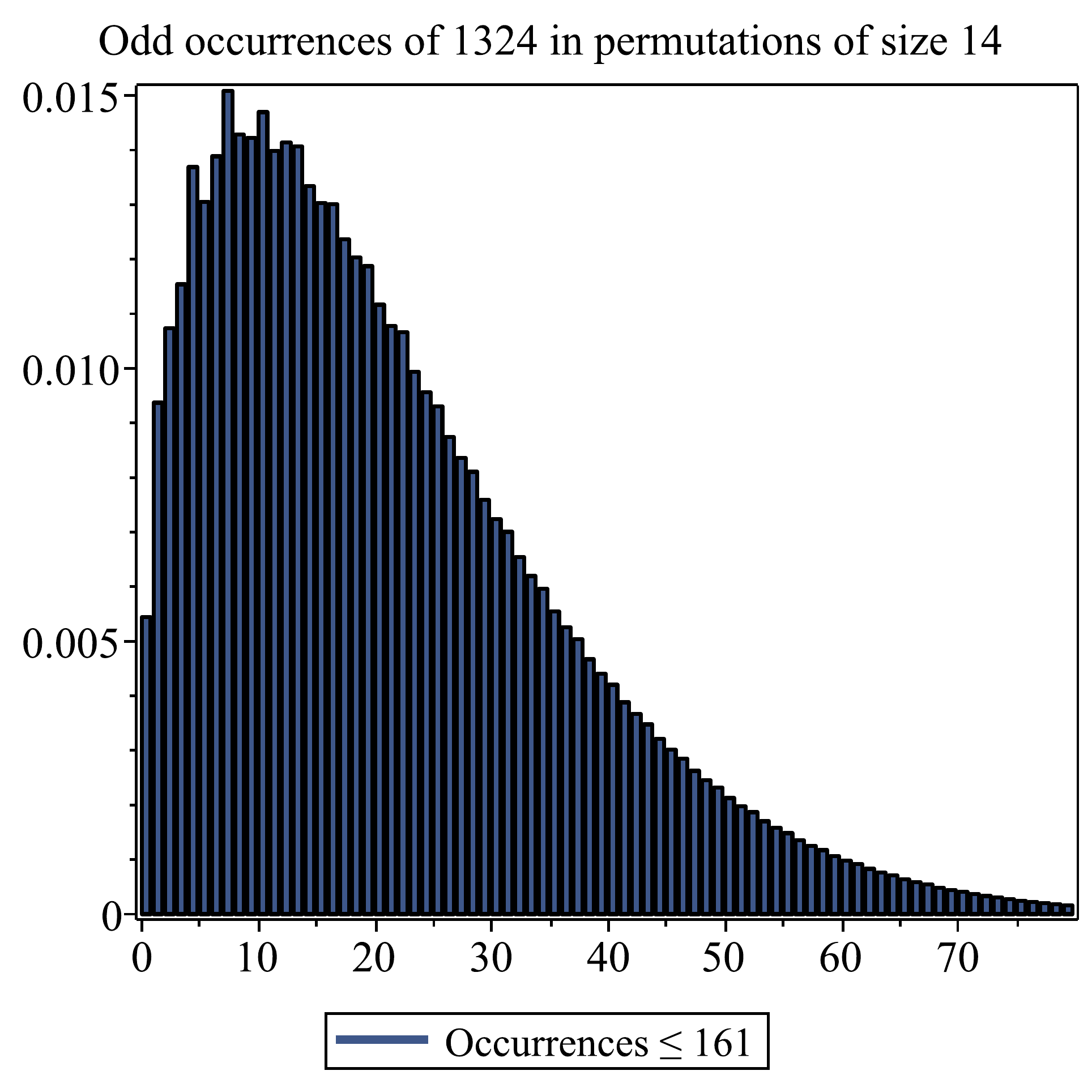}}
\caption{Histogram of 1324 odd occurrences.} 
\label{fig:1324o}
\end{minipage}
\hspace{0.05\textwidth}
\begin{minipage}[t]{0.45\textwidth} 
\centerline {\includegraphics[width=\textwidth]{ 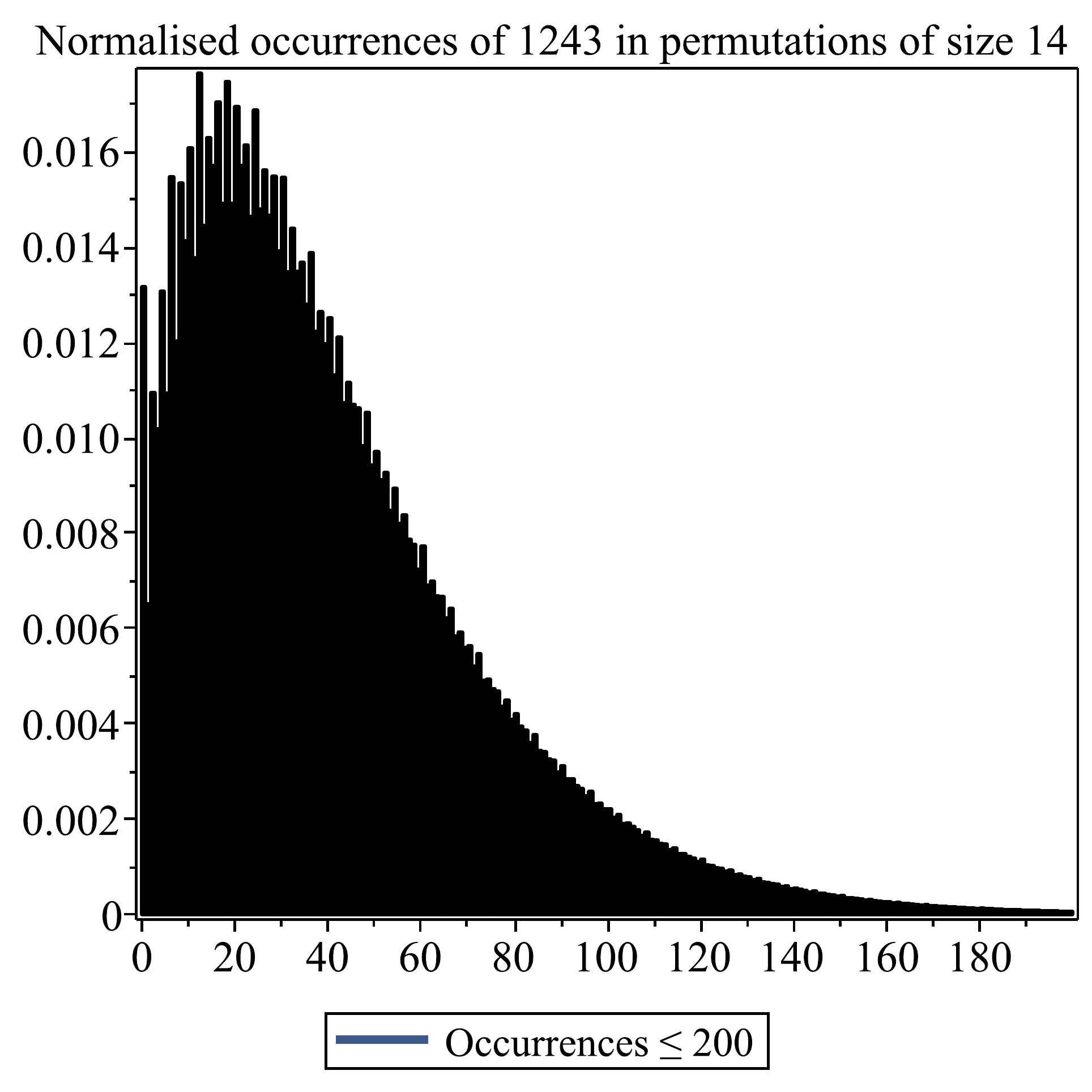}}
\caption{Histogram of 1243 occurrences.} 
\label{fig:1243}
\end{minipage}
\end{figure}

For Wilf classes VI and VII there are no pre-existing data for $r > 0.$ For Wilf class VI we find:
For $r=1,$
\begin{align*}
  & 0,0,0,1,10,77,548,3799,26165,180512,1251832,8738589,61427007,\\
  &434771094,3097485378,22203860315,160077190385.
  \end{align*}
And for $r=2,$
\begin{align*}  
 &0, 0, 0, 0, 6, 69, 598, 4686, 35148, 258390, 1882813, \\
  &13677083, 99350385, 722871146, 5272996671,
\end{align*}
and for Wilf class VII for $r=1,$
\begin{align*}
& 0, 0, 0, 1, 9, 62, 402, 2593, 16921, 112196, 755920, 5168174, 35796046, \\
  &250765372, 1774228404, 12662584870, 91064282806.
  \end{align*}
  And for $r=2$
  \begin{align*}
 &0, 0, 0, 0, 8, 82, 612, 4187, 28065, 188514, 1278590, \\
  &8774123, 60914835, 427488844, 3029373540.
\end{align*}
We have shown in Figs. \ref{fig:1234}, \ref{fig:1342}, \ref{fig:1324}, \ref{fig:1243} the normalised histograms for 4 of the 7 Wilf classes.
The only reason we don't show the histograms for the remaining 3 classes is that they provide no further illumination. Unlike the case for patterns of length 3, there appears to be an odd-even parity effect in every case, so we have also shown the odd- and even-subsequences for these four classes. Also unlike the situation for patterns of length 3, the histograms are quite skewed, compared to the comparatively symmetric shape (apart from the long, skinny tail) in the length-3 case.

As these distributions are provably asymptotically normal \cite{JNZ13}, the observed parity effects and asymmetry must be a manifestation of permutations that are too short. That is to say, we are far from the asymptotic regime.

Having observed the two different $r$-dependencies in the length-3 situation, we now investigate the corresponding asymptotics for the seven length-4 Wilf classes.

\subsection{Class I}
As remarked above, for $r=1$ we have 200 terms, and for $r=2$ we have 25 terms. As we only have 17 terms for some of the classes that haven't previously been studied, we will first confine our analysis to 17-term sequences, and then see if the additional terms change our conclusions. We will show our analysis in some detail for this case, and then apply exactly the same methods in all subsequent cases, and just quote the results.

As our analysis will be based on the ratio method and its extensions, as described in the appendix, it would be useful to have as many terms as possible, even if these terms are only approximate. 

The method of series extension \cite{G16} allows us to predict additional terms by constructing  holonomic differential equations from the known coefficients -- known as the method of differential approximants \cite{GJ72} -- and using these to predict further terms. This method is also briefly described in the appendix. By this method the accuracy of the predicted terms is assessed by constructing the variance of the estimates of each predicted coefficient, which we are able to do as we construct many distinct differential equations by varying the degrees of the polynomials multiplying the various derivatives. We cut off the number of predicted terms when the estimated error grows larger than some desired maximum.

In this case we were able to predict 12 further coefficients, as shown in Table \ref{tab1} below, where for comparison we have also shown the exact coefficients. It can be seen that the first predicted coefficient is in error by 1 part in the 9th significant digit, while the last predicted coefficient is in error by 6 parts in the 6th significant digit. These approximate coefficients are still useful for our analysis, as we now show.

\begin{figure}[ht!] 
\begin{minipage}[t]{0.45\textwidth} 
\centerline{\includegraphics[width=\textwidth]{ 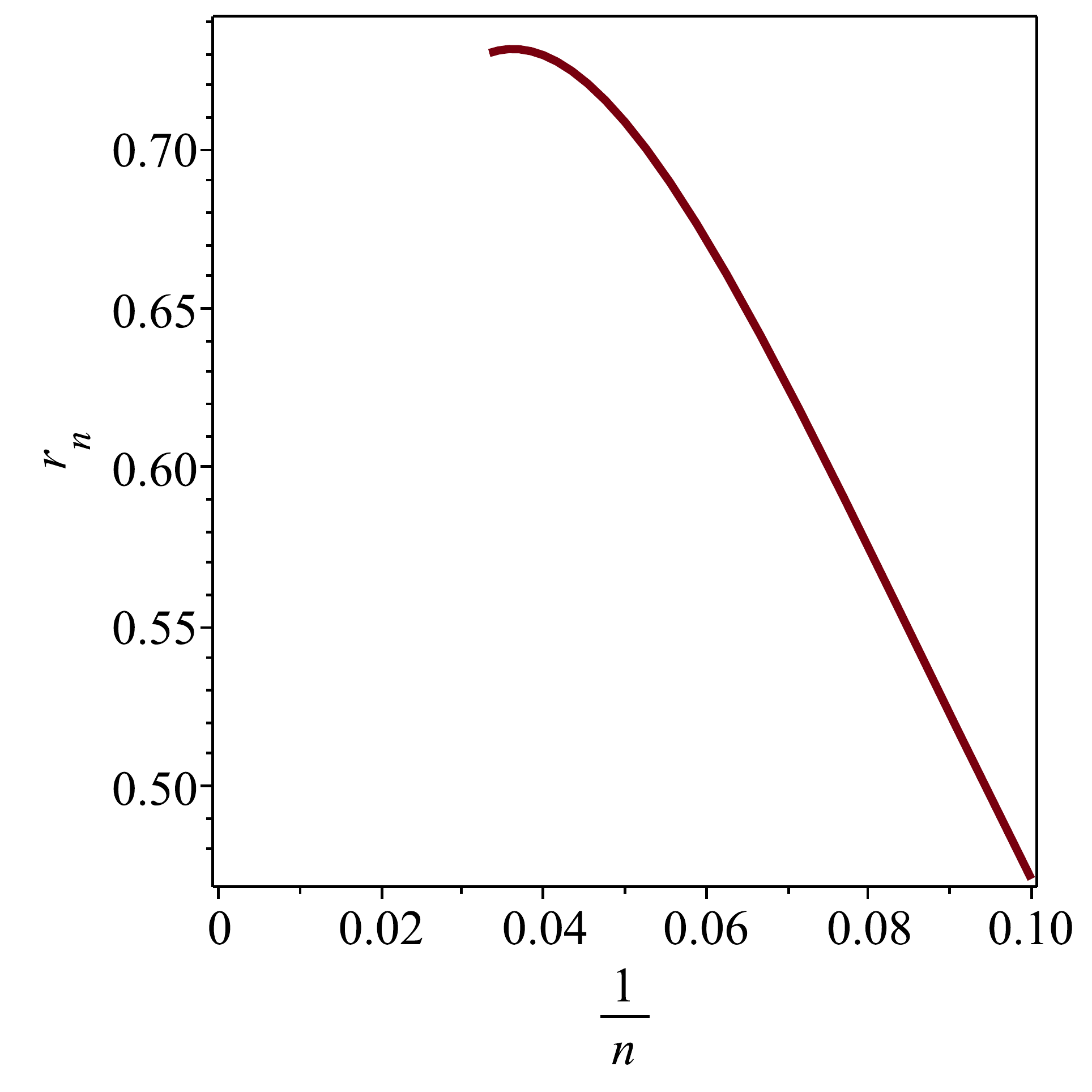}}
\caption{Ratio of $\psi_1(n)/\psi_0(n)$ plotted against $1/n$ for 1234 patterns.} 
\label{fig:1234r1}
\end{minipage}
\hspace{0.05\textwidth}
\begin{minipage}[t]{0.45\textwidth} 
\centerline {\includegraphics[width=\textwidth]{ 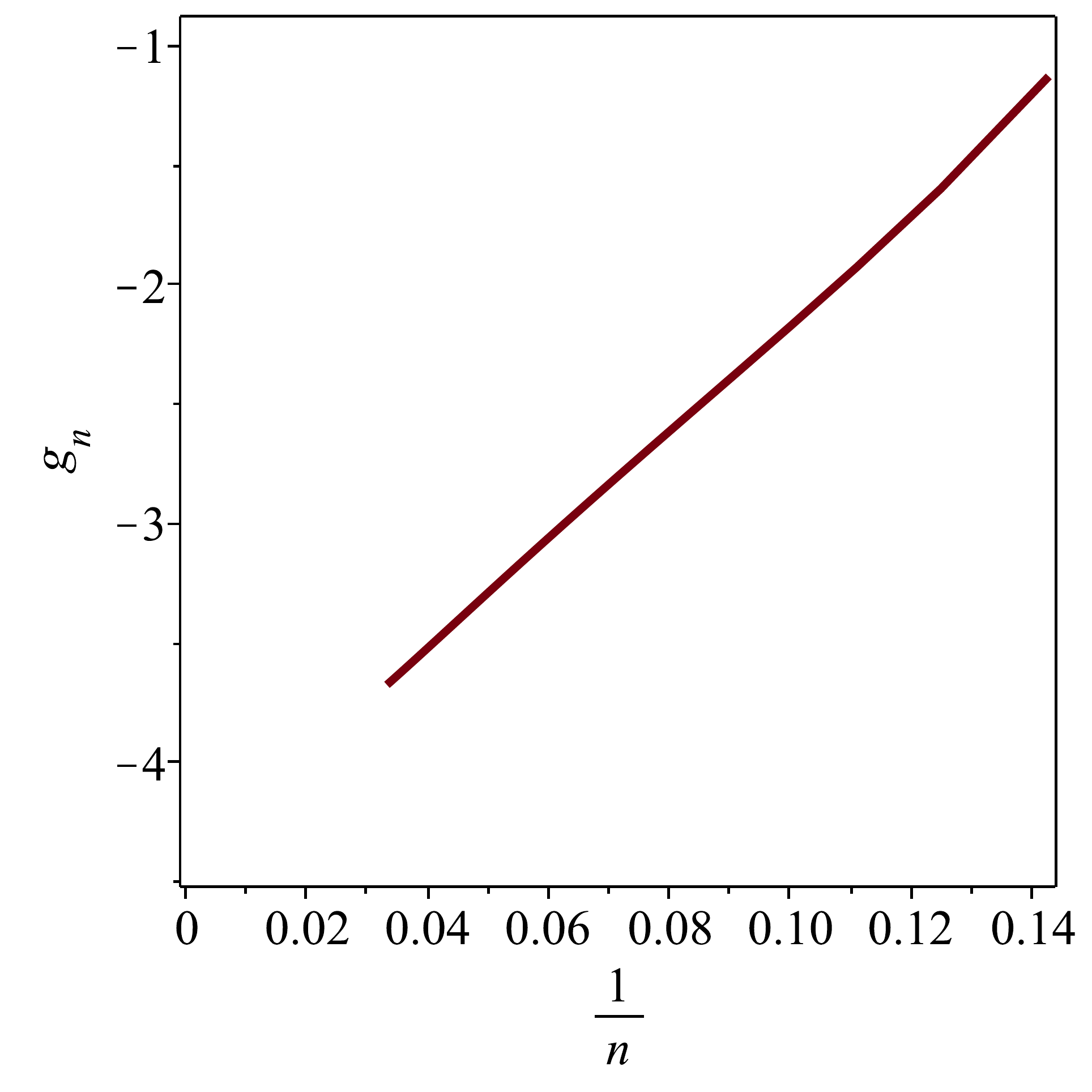}}
\caption{Estimators $g_n$ of exponent $g$ plotted against $1/n.$} 
\label{fig:1234g1}
\end{minipage}
\end{figure}

Recall that for the pattern 123 the asymptotics did not change with the number of occurrences, only the amplitude did. We know that the exponential growth does not change, so the only possible change is in the power-law exponent. Recall that $\psi_0(n)\sim\frac{81\sqrt{3}\cdot 9^n}{16\pi\cdot n^4},$ see A005802 \cite{OEIS}. So $\psi_1(n)\sim\frac{C_1\cdot 9^n}{16\pi\cdot n^g},$ where $C_1$ and the exponent $g$ are unknown. If $g=4,$ then the ratio $$\frac{\psi_0(n)}{\psi_1(n) }\sim \frac{ 81\sqrt{3}}{C_1\cdot n^{4-g}}$$ should tend to a non-zero constant $81\sqrt{3}/C_1$ as $n$ increases. Otherwise the ratio should diverge or vanish, according as $g$ is less than or greater than 4. It can be seen from Fig. \ref{fig:1234r1} that the ratio is certainly not diverging, but rather appears to be going to a constant.
We can also use the ratio method to estimate the value of $g$ directly, as the ratios $r_n= \psi_1(n)/\psi_1(n-1) \sim  9(1-g/n+o(1/n)). $ We estimate the exponent $g$ from the estimators $g_n =(r_n/9-1)n \sim -g(1+o(1)).$

\begin{figure}[ht!] 
\begin{minipage}[t]{0.45\textwidth} 
\centerline{\includegraphics[width=\textwidth]{ 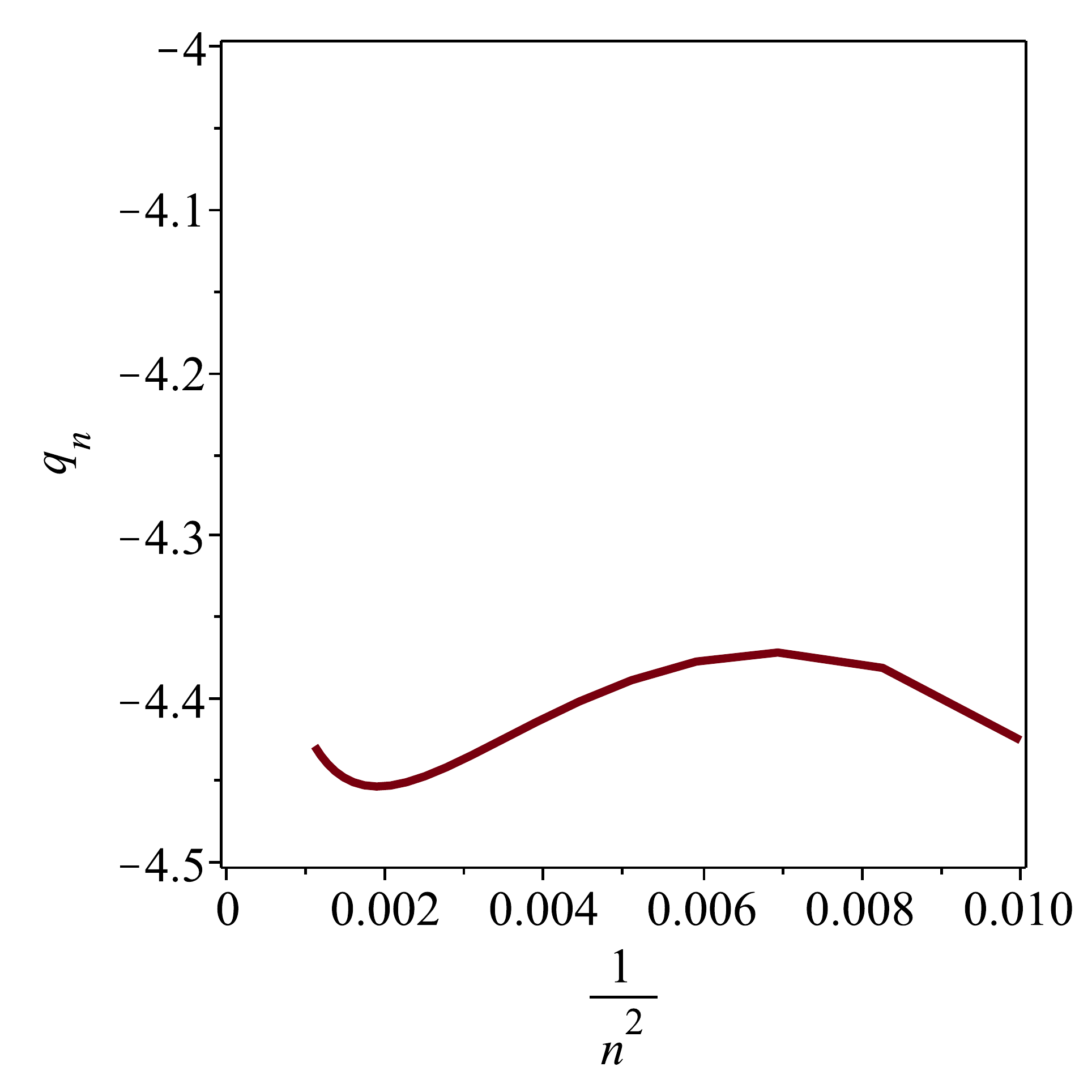}}
\caption{Quadratic estimators $q_n$ of exponent $g$ plotted against $1/n^2,$ using 30 terms.} 
\label{fig:1234g2}
\end{minipage}
\hspace{0.05\textwidth}
\begin{minipage}[t]{0.45\textwidth} 
\centerline {\includegraphics[width=\textwidth]{ 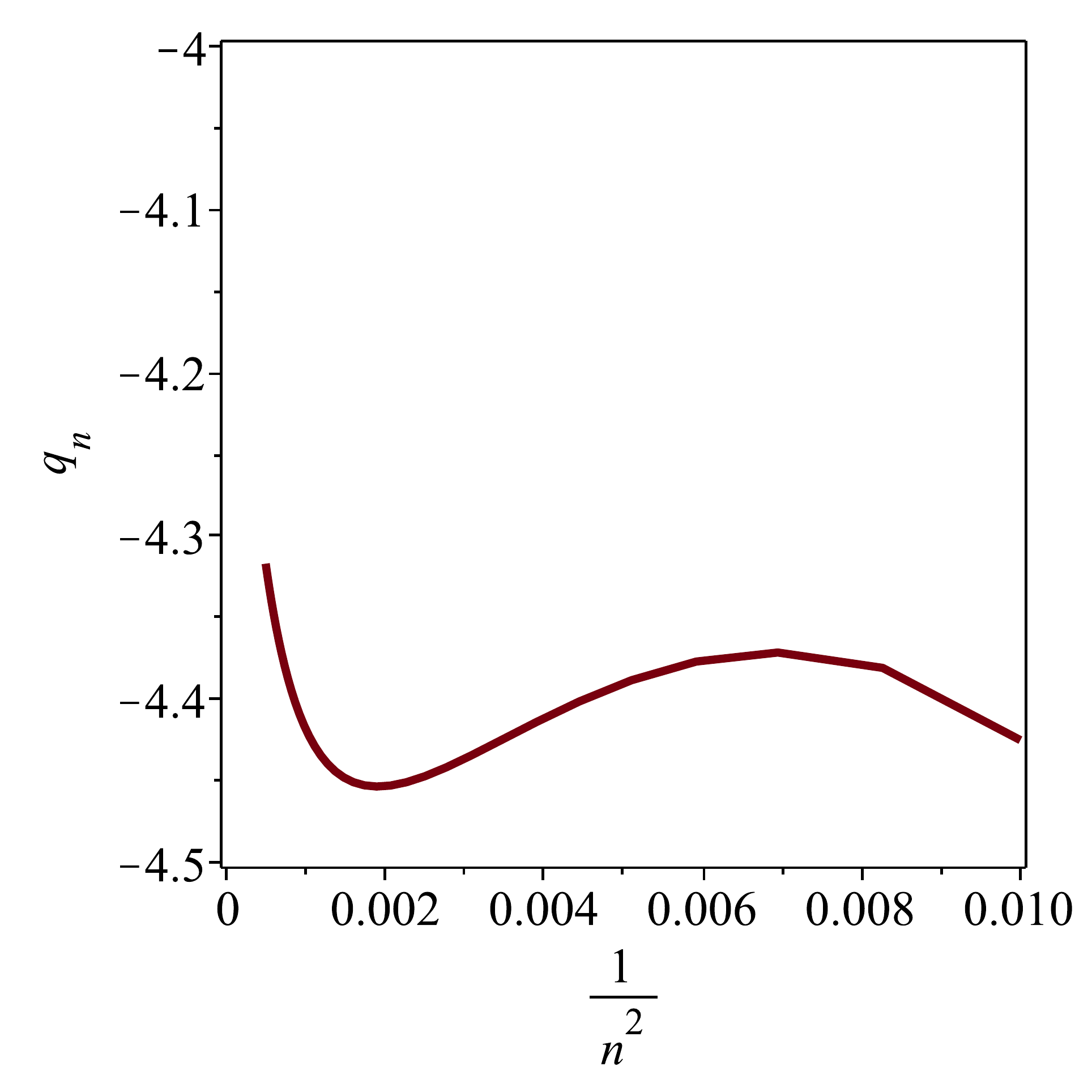}}
\caption{Quadratic estimators $q_n$ of exponent $g$ plotted against $1/n^2,$ using 45 terms.} 
\label{fig:1234g3}
\end{minipage}
\end{figure}

In Fig. \ref{fig:1234g1} we plot the estimators $g_n$ against $1/n.$ These are seen to be going to a value close to -4, though possibly a little lower. If the singularity is a pure power-law singularity, the sub-dominant term $o(1)$ is in fact $O(1/n).$ We can eliminate this term by forming the quadratic estimators $q_n \equiv ng_n - (n-1)g_{n-1} \sim g(1-o(1/n)).$
In Fig. \ref{fig:1234g2} we show these quadratic estimators plotted against $1/n^2,$ and it can be seen that there is an upturn in the estimators away from -4.5 and toward the value -4. To show that this is not spurious, we show in Fig. \ref{fig:1234g3} the same plot, but this time using 45 terms, which we have in this case, and which clearly shows that the limiting value -4 is quite compelling.

\begin{figure}[ht!] 
\begin{minipage}[t]{0.45\textwidth} 
\centerline{\includegraphics[width=\textwidth]{ 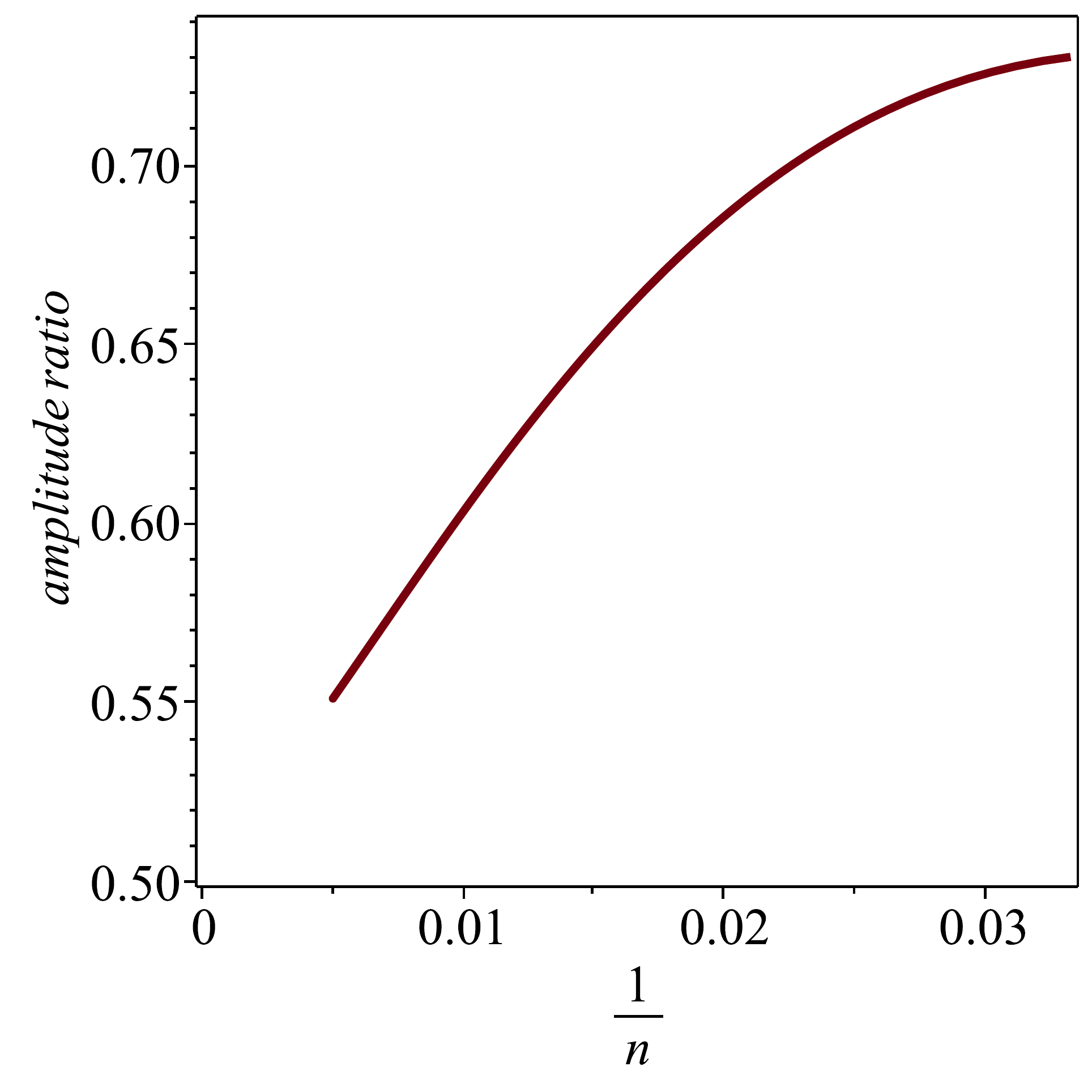}}
\caption{Amplitude ratios $\psi_1(n)/\psi_0(n)$ plotted against $1/n,$ using 200 terms.} 
\label{fig:1234r2}
\end{minipage}
\hspace{0.05\textwidth}
\begin{minipage}[t]{0.45\textwidth} 
\centerline {\includegraphics[width=\textwidth]{ 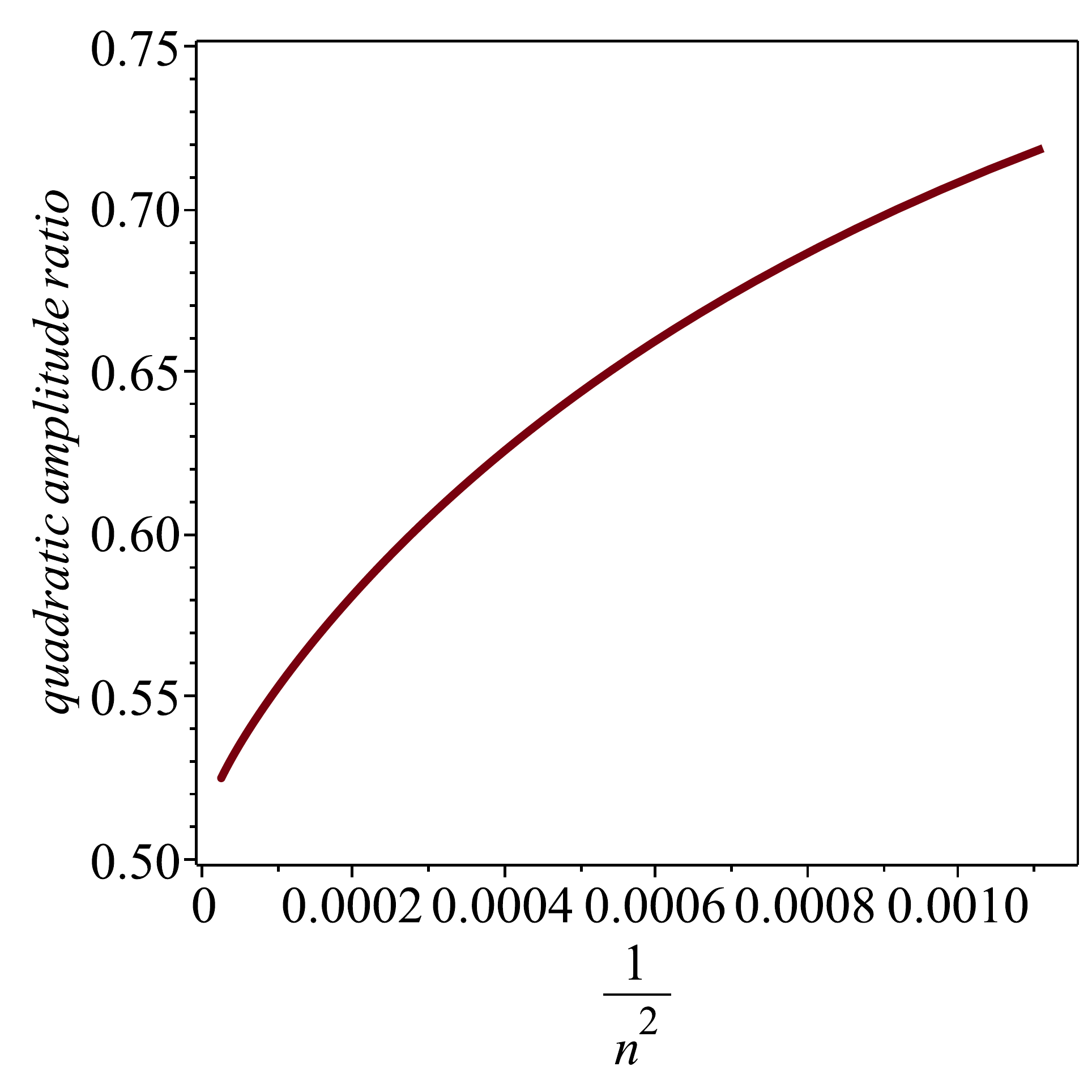}}
\caption{Quadratic amplitude ratios plotted against $1/n^2,$ using 200 terms.} 
\label{fig:1234r3}
\end{minipage}
\end{figure}

It remains to estimate the value of the amplitude $C_1.$ This is equivalent to estimating the limit of the plot shown in Fig. \ref{fig:1234r1}, which we can't do with the 29 terms we have at our disposal. However if we now use all 200 known terms, then this ratio is clearly going to a limit around $0.5,$ as can be seen from Fig. \ref{fig:1234r2}. And in Fig. \ref{fig:1234r3} we show quadratic estimators of the same ratio, plotted against $1/n^2.$ It is clear that the limit is plausibly 0.5, which implies that 
$$\psi_1(n)\sim \frac{81\sqrt{3}\cdot 9^n}{32\pi\cdot n^4}.$$ Recall from the previous section that for the pattern 123 we have 
$\psi_1(n)/\psi_0(n)=3.$ For the pattern 1234 we conjecture that $\psi_1(n)/\psi_0(n)=1/2.$

For the 2-occurrence sequence, as noted above, we have 25 terms available. By the method of series extension we are able to obtain a further 40 approximate terms. We have analysed that sequence exactly as described above for the 1-occurrence sequence, and we find exactly the same asymptotics, except that the ratio $\psi_2(n)/\psi_0(n)\approx 1.$ So we conjecture that 
$$\psi_2(n)\sim \psi_0(n) \sim \frac{81\sqrt{3}\cdot 9^n}{16\pi\cdot n^4}.$$

It seems reasonable to conjecture that $\psi_r(n)=\frac{C_r\cdot 9^n}{n^4}$ for all $r,$ just as is the case for the shorter pattern 123.
That is to say, as $r$ changes, only the amplitude $C_r$ changes. The exponential growth remains (provably) the same, and the sub-dominant power law term is also (conjecturally) unchanged. The pre-multiplicative amplitude is also of course just conjectured, rather than proved.

\begin{table}[htp]
\caption{Predicted and exact coefficients of order 18 to 29}
\vspace{2mm}
\centering
\begin{tabular}{|c|c|c|}
\hline
Order & Predicted& Exact\\
\hline
18&2.05621895e+12&2056218941678\\
19&1.53582968e+13&  15358296210724\\
20&1.15469577e+14&  115469557503753\\
21&8.73561640e+14&  873561194459596\\
22&6.6477690e+15&  6647760790457218 \\
23&5.0871661e+16&50871527629923754\\
24&3.9134709e+17&  391345137795371013\\
25&3.0255951e+18&  3025568471613091692\\
26&2.3502067e+19&  23501724670464335914\\
27&1.8337471e+20&  183370520135071994536\\
28&1.4368444e+21&  1436795093911521996331\\
29&1.1303753e+22&11303188383039278887124\\
\hline
\end{tabular}
\label{tab1}
\end{table}%

\begin{figure}[ht!] 
\begin{minipage}[t]{0.45\textwidth} 
\centerline{\includegraphics[width=\textwidth]{ 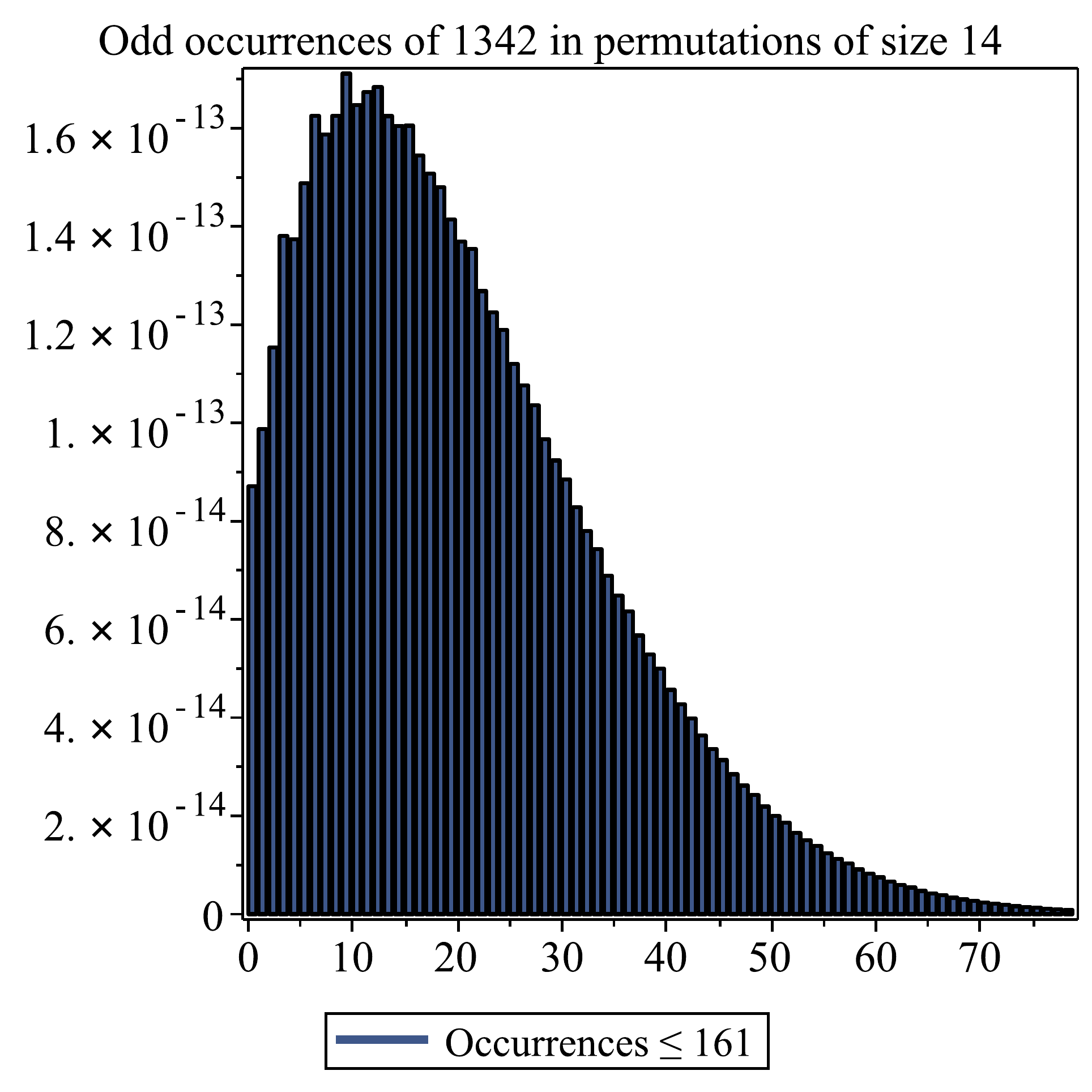}}
\caption{Histogram of 1342 odd occurrences.} 
\label{fig:1342o}
\end{minipage}
\hspace{0.05\textwidth}
\begin{minipage}[t]{0.45\textwidth} 
\centerline {\includegraphics[width=\textwidth]{ 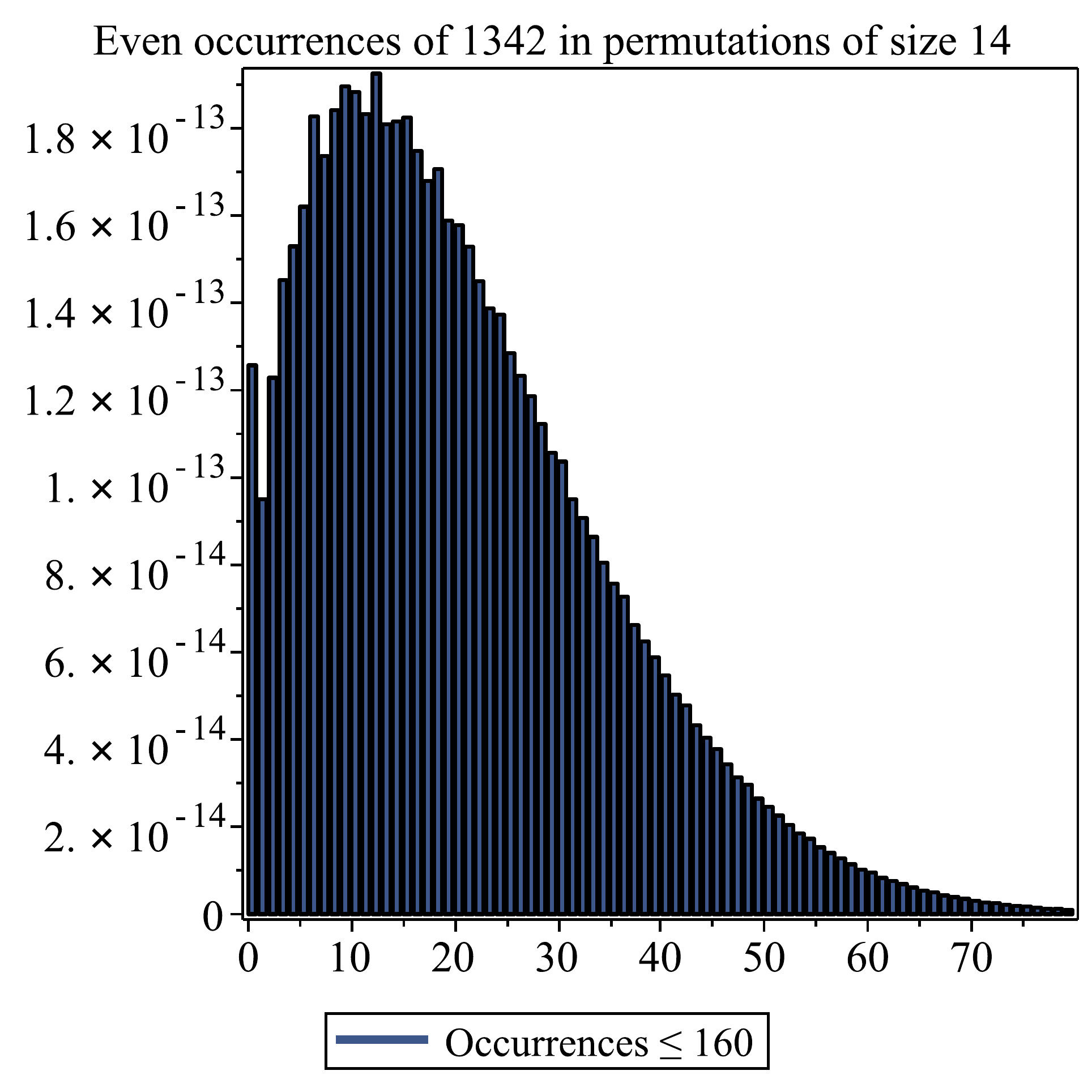}}
\caption{Histogram of 1342 even occurrences.} 
\label{fig:1342e}
\end{minipage}
\end{figure}

\begin{figure}[ht!] 
\begin{minipage}[t]{0.45\textwidth} 
\centerline{\includegraphics[width=\textwidth]{ 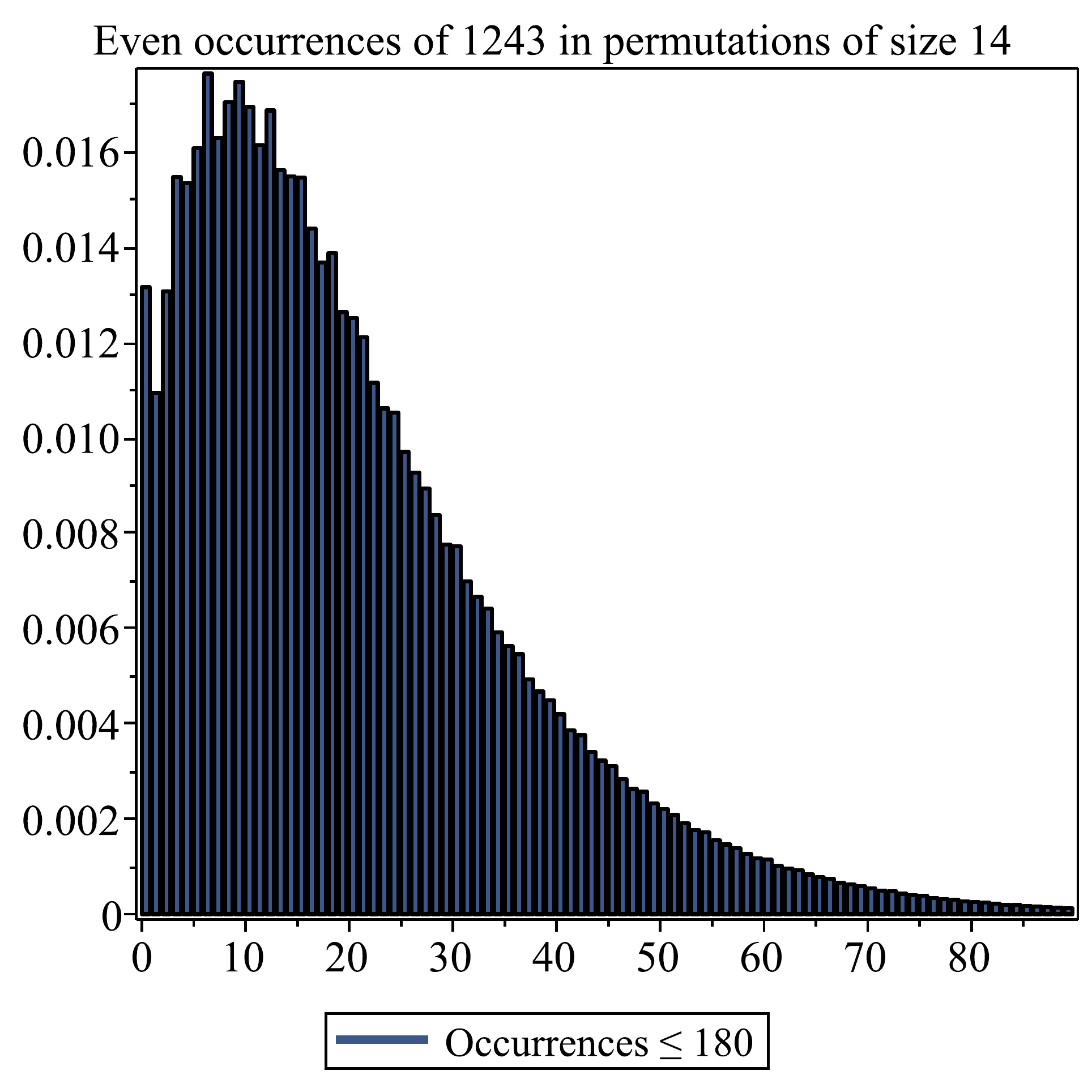}}
\caption{Histogram of 1243 even occurrences.} 
\label{fig:1243e}
\end{minipage}
\hspace{0.05\textwidth}
\begin{minipage}[t]{0.45\textwidth} 
\centerline {\includegraphics[width=\textwidth]{ 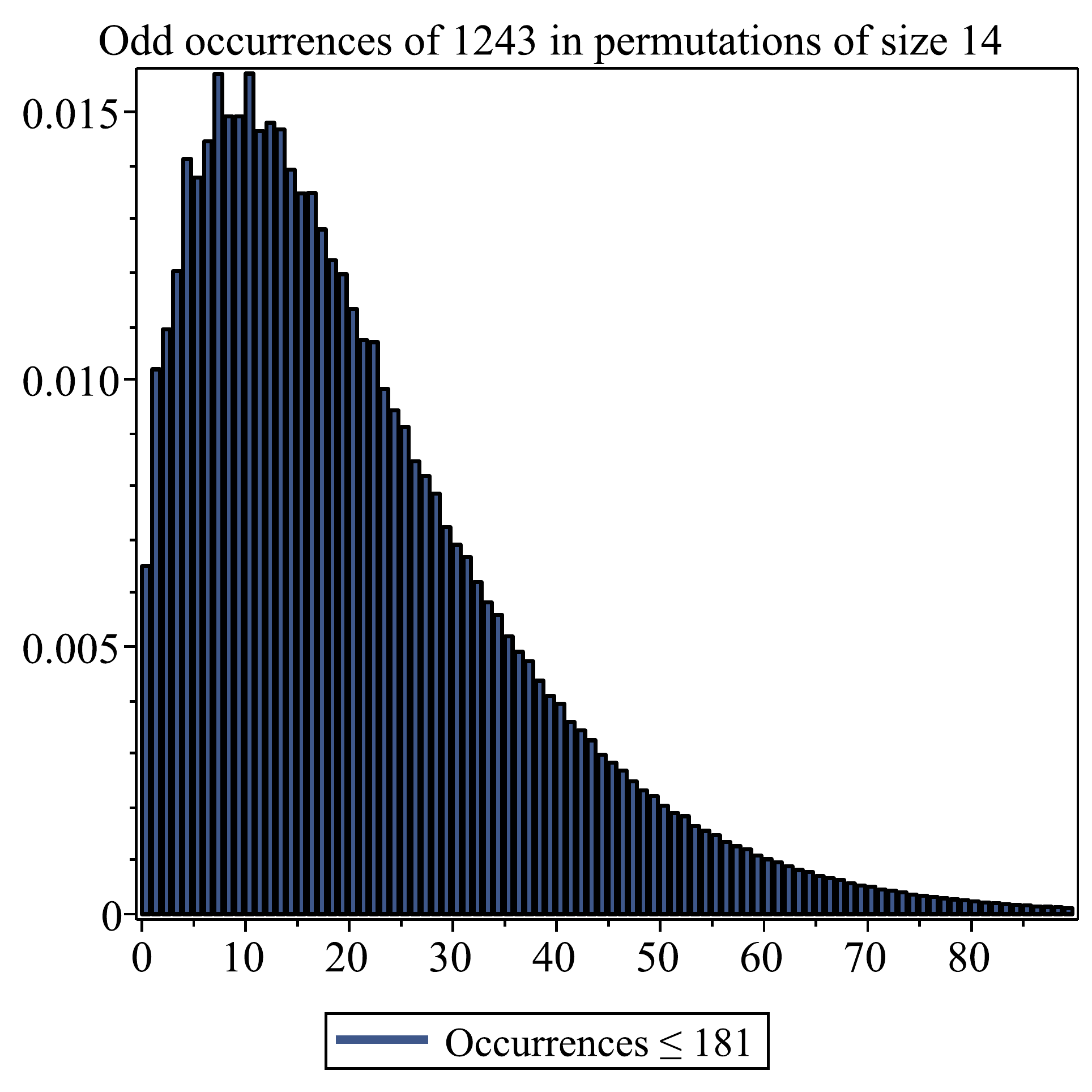}}
\caption{Histogram of 1243 odd occurrences.} 
\label{fig:1243o}
\end{minipage}
\end{figure}

\newpage
\subsection{Class II}
We analysed OEIS sequence A224179 for 1-occurrences of the pattern 1243. 25 terms are given, and we were able to find a further 50 approximate terms with sufficient apparent accuracy for the ratio method to be used.
We analysed this sequence in an identical manner to our analysis of 1234 1-occurrences, discussed in the previous section. 
The ratio $\psi(1243)_1(n)/\psi(1234)_0(n)$ is going to an apparent limit of $0.93125 \pm 0.00015.$ If this were a simple fraction, $149/160 =0.93125$ would seem to be a likely candidate.

Direct analysis of the coefficients of the sequence $\psi(1243)_1(n)$ for the power-law exponent also verified the O$(n^{-4})$ behaviour that follows from the fact that the ratio we have just estimated is non-zero and finite.

We conclude that $\psi_1(n)\sim \frac{C_1\cdot 9^n} {n^4},$ where $C_1\approx \frac{0.93125 \cdot 81\sqrt{3}}{16\pi},$ and perhaps $C_1 = 149\cdot 3^4\sqrt{3}/(2^9\cdot 5\pi)$ exactly.

We similarly analysed the OEIS sequence A224181 for 2-occurrences of the pattern 1243. 24 terms are given, and we were able to find a further 40 approximate terms with sufficient apparent accuracy for the ratio method to be used. We analysed this sequence in an identical manner to our analysis of 1243 1-occurrences, just discussed.
The ratio $\psi(1243)_2(n)/\psi(1234)_0(n)$ is going to an apparent limit of $2.93 \pm 0.01.$ It is tempting to guess that this is just 2 greater than the corresponding ratio for the 1-occurrence case. If that were so, then a simple fraction, $469/160 =2.93125$ would seem to be a likely candidate.

We conclude that
$\psi_2(n)\sim \frac{C_2\cdot 9^n} {n^4},$ where $C_2\approx \frac{2.93 \cdot 81\sqrt{3}}{16\pi},$ and perhaps $C_2 = 469\cdot 3^4\sqrt{3}/(2^9\cdot 5\pi)$ exactly.

As for the previous class, it seems reasonable to conjecture that $\psi_r(n)=\frac{C_r\cdot 9^n}{n^4}$ for all $r,$ just as is the case for the shorter pattern 123.
That is to say, as $r$ changes, only the amplitude $C_r$ changes. The exponential growth remains (provably) the same, and the sub-dominant power law term is also (conjecturally) unchanged.

\subsection{Class III}
We analysed the OEIS sequence A224182 for 1-occurrences of the pattern 1432. 18 terms are given, and we were able to find a further 26 approximate terms with sufficient apparent accuracy for the ratio method to be used.
We analysed this sequence in a similar manner to our analysis of 1234 1-occurrences, discussed above. 
However the ratio $\psi(1432)_1(n)/\psi(1234)_0(n)$ appears to be diverging linearly. Accordingly, we studied the sequence $\psi(1432)_1(n)/(n\cdot \psi(1234)_0(n)).$ That ratio is going to an apparent limit of $0.0375 \pm 0.0005.$ It is tempting to guess that this is just $3/80.$ 

We conclude that
$\psi_1(n)\sim \frac{C_1\cdot 9^n} {n^3},$ where $C_1\approx \frac{0.0375 \cdot 81\sqrt{3}}{16\pi},$ and perhaps\\ $C_1 =  3^5\sqrt{3}/(2^8\cdot 5\pi)$ exactly.

We next analysed the sequence for 2-occurrences of the pattern 1432. We have generated 15 terms, and we were able to find a further 9 approximate terms with sufficient apparent accuracy for the ratio method to be used. We analysed this sequence in an identical manner to our analysis of 1432 1-occurrences, just discussed.
The ratio $\psi(1432)_2(n)/(n^2\cdot \psi(1234)_0(n))$ is going to an apparent limit of $0.00075 \pm 0.00005.$ This is too imprecise to hazard a guess at the exact value.

We conclude that
$$\psi_2(n)\sim \frac{C_2\cdot 9^n} {n^2},$$ where $C_2\approx \frac{0.00075 \cdot 81\sqrt{3}}{16\pi}.$ 

We see that the behaviour of the $r$-occurrences is qualitatively similar to that of 132-avoiders, in that the power-law exponent increases by 1 as $r$ increases by 1. Accordingly, we conjecture that $$\psi_r(n)\sim \frac{C_r\cdot 9^n} {n^{4-r}},$$ where $C_0$ is known, $C_1$ is conjectured, and $C_2$ is estimated.

\subsection{Class IV}
There is no pre-existing data for this class, of which 2143 is a representative member. We have generated 18 coefficients for $\psi_1(n)$ and 16 coefficients for $\psi_2(n).$ We were able to find a further 20 approximate terms for $\psi_1(n)$  and 15 approximate terms for $\psi_2(n)$ with sufficient apparent accuracy for the ratio method to be used.

We analysed these sequences in a similar manner to our analysis of other sequences, discussed above. 
The ratio $\psi(2143)_1(n)/\psi(1234)_0(n)$ again appears to be diverging linearly. Accordingly, we studied the sequence $\psi(2143)_1(n)/(n\cdot \psi(1234)_0(n)).$ That ratio is going to an apparent limit of $0.0421 \pm 0.0005.$ 

We conclude that $\psi_1(n)\sim \frac{C_1\cdot 9^n} {n^3},$ where $C_1\approx \frac{0.0421 \cdot 81\sqrt{3}}{16\pi}.$ 

We next analysed the sequence for 2-occurrences of the pattern 2143. We analysed this sequence in an identical manner to our analysis of 2143 1-occurrences, just discussed.
The ratio $\psi(2143)_2(n)/(n^2\cdot \psi(1234)_0(n))$ is going to an apparent limit of $0.00090 \pm 0.00005.$ 

We conclude that
$\psi_2(n)\sim \frac{C_2\cdot 9^n} {n^2},$ where $C_2\approx \frac{0.0009 \cdot 81\sqrt{3}}{16\pi}.$ 

We see that for this class too, the behaviour of the $r$-occurrences is qualitatively similar to that of 132-avoiders, in that the power-law exponent increases by 1 as $r$ increases by 1. Accordingly, we conjecture that $$\psi_r(n)\sim \frac{C_r\cdot 9^n} {n^{4-r}},$$ where $C_0$ is known, and $C_1$ and $C_2$ are estimated.
\subsection{Class V}
For this class, of which 1324 is a representative member, 17 terms for $\psi_1(n)$ are given as OEIS sequence A224182, but there is no pre-existing data for $\psi_2(n).$ We have given 15 terms for that sequence. We have also used a further 15 and 8 approximate terms in our analysis, for 1- and 2-occurrences respectively.

Recall that for this class, uniquely, the result for 0-occurrences is not known rigorously. Numerical work by Conway et al. \cite{CGZ18} suggests that $\psi_0(n) \sim C_0 \cdot \mu^n \cdot \mu_1^{\sqrt{n}} \cdot n^g,$ where $\mu \approx 11.598$ (and possibly $9+3\sqrt{3}/2$ exactly), $\mu_1 \approx 0.0400,$ and $g \approx -1.1.$ 

Recall that $\psi_r(n)$ grows at the same rate as $\psi_0(n),$ so we studied the ratios $$\psi_r(n)/\psi_0(n)$$ for $r=1$ and $r=2.$

The ratio $\psi(1324)_1(n)/\psi(1324)_0(n)$ appears to be diverging linearly. Accordingly, we studied the sequence $\psi(1324)_1(n)/(n\cdot \psi(1324)_0(n)).$ That ratio is going to an apparent limit of $0.013 \pm 0.001.$ 

We conclude that
$$\psi_1(n)\sim C_1\cdot n\cdot \psi_0(n),$$ where $C_1 = 0.013 \pm 0.001.$ 

For 2-occurrences, we studied the sequence $\psi(1324)_2(n)/(n^2\cdot \psi(1324)_0(n)).$ That ratio is going to an apparent limit of $0.0010 \pm 0.0001.$ So we conclude that
$\psi_2(n)\sim C_2\cdot n^2\cdot \psi_0(n),$ where $C_2 = 0.001 \pm 0.0001.$ 

We see that for this class too, the behaviour of the $r$-occurrences is qualitatively similar to that of 132-avoiders, in that the power-law exponent apparently increases by 1 as $r$ increases by 1. Accordingly, we conjecture that $$\psi_r(n)\sim C_r\cdot n^{r}\cdot \psi_0(n).$$

\subsection{Class VI}
For this class and class VII we have $\psi_0(n) \sim \frac{64\cdot 8^n}{243\sqrt{\pi}\cdot n^{5/2}}.$ There are no pre-existing sequences for any $r > 0.$ For both classes we have generated 17 terms for 1-occurrences and 15 terms for 2-occurrences. We have also generated 30 approximate terms for 1-occurrences and 8 approximate terms for 2-occurrences, for both classes.

The ratio $\psi(1342)_1(n)/\psi(1342)_0(n)$ again appears to be diverging linearly. Accordingly, we studied the sequence $\psi(1342)_1(n)/(n\cdot \psi(1342)_0(n)).$ That ratio is going to an apparent limit of $0.0206 \pm 0.0005.$ 

We conclude that
$$\psi_1(n)\sim \frac{C_1\cdot 8^n} {n^{3/2}},$$ where $C_1=0.0206 \pm 0.0005.$ 

We analysed the sequence for 2-occurrences in an identical manner to our analysis of 1342 1-occurrences, just discussed.
The ratio $\psi(1342)_2(n)/(n^2\cdot \psi(1342)_0(n))$ is going to an apparent limit of $0.0005 \pm 0.0001.$ So we conclude that
$$\psi_2(n)\sim \frac{C_2\cdot 8^n} {n^{1/2}},$$ where $C_2=0.0005 \pm 0.0001.$ 

\subsection{Class VII}
The ratio $\psi(2413)_1(n)/\psi(2413)_0(n)$ again appears to be diverging linearly. We again studied the sequence $\psi(2413)_1(n)/(n\cdot \psi(2413)_0(n)).$ That sequence is going to an apparent limit of $0.01463 \pm 0.00005.$ 

We conclude that
$$\psi_1(n)\sim \frac{C_1\cdot 8^n} {n^{3/2}},$$ where $C_1=0.01463 \pm 0.00005.$ 

We analysed the sequence for 2-occurrences in an identical manner to our analysis of 2413 1-occurrences, just discussed.
The ratio $\psi(2413)_2(n)/(n^2\cdot \psi(2413)_0(n))$ is going to an apparent limit of $0.00010 \pm 0.00002.$ So we conclude that
$$\psi_2(n)\sim \frac{C_2\cdot 8^n} {n^{1/2}},$$ where $C_2=0.00010 \pm 0.00002.$ 

We see that for both class VI and class VII the behaviour of the $r$-occurrences is qualitatively similar to that of 132-avoiders, in that the power-law exponent apparently increases by 1 as $r$ increases by 1. Accordingly, we conjecture that $$\psi_r(n)\sim C_r\cdot n^{r}\cdot \psi_0(n)$$ for both classes.

\subsection{Class summary}
Just as for patterns of length three, we see that the behaviour of $r$-occurrences of patterns of length four falls into two groups.

The first group consists of classes I and II. For those two classes, the asymptotic behaviour is unchanged as $r$ increases, and only the pre-multiplicative constant $C_r$ changes with $r$. That is to say, $\psi_r(n)\sim C_r\cdot \psi_0(n).$

For the second group, consisting of classes III-VII, the subdominant power-law behaviour changes, increasing by 1 with each unit increase in $r.$ It is also the case that the pre-multiplicative constant $C_r$ changes with $r$. That is to say, $\psi_r(n)\sim C_r\cdot n^{r}\cdot \psi_0(n).$

\section{Conclusion}
We have given details of a new algorithm for enumerating multisets that is likely to be of broad applicability. As a first application, we have studied the behaviour of occurrences of patterns in permutations. For patterns of length 3 we give the general solution for the two Wilf classes, and have given the first systematic study of patterns of length 4, based on extensive enumerations obtained from our algorithm, combined with some longer enumerations obtained by others in isolated special cases. We find that there are seven Wilf classes, and have given conjectured asymptotic formulae for the behaviour of each of these classes.

\section{Acknowledgements} We are particularly grateful to Yuma Inoue who modified his program for counting PAPs to counting occurrences of patterns, and using his program we were able to
provide verification of the results produced by our program, and in a couple of specific cases obtain an additional term in the sequence.  We would also like to thank Vince Vatter and Doron Zeilberger for helpful comments on the manuscript.

\end{document}